\documentclass[11pt]{article}

\usepackage[margin=1.1in]{geometry}
\usepackage{amsmath,amssymb,amsthm,mathtools}
\usepackage{enumitem}
\usepackage{multirow}
\usepackage{authblk}
\usepackage[colorlinks=true,linkcolor=blue,citecolor=blue,urlcolor=blue]{hyperref}

\newtheorem{theorem}{Theorem}[section]
\newtheorem{lemma}[theorem]{Lemma}

\newtheorem{corollary}[theorem]{Corollary}
\theoremstyle{definition}
\newtheorem{definition}[theorem]{Definition}
\newtheorem{conjecture}[theorem]{Conjecture}
\theoremstyle{remark}
\newtheorem{remark}[theorem]{Remark}

\newcommand{\R}{\mathbb R}
\newcommand{\Z}{\mathbb Z}
\newcommand{\norm}[1]{\left\lVert#1\right\rVert}
\newcommand{\abs}[1]{\left\lvert#1\right\rvert}
\newcommand{\J}{\mathcal J}
\newcommand{\W}{\mathcal W}
\newcommand{\Pk}{\mathcal P_{k}}
\newcommand{\Hsp}{H^{s}}
\newcommand{\Hm}{H^{-s}}
\newcommand{\Bball}{B^{s}}
\newcommand{\dn}{d_{n}}
\newcommand{\Lam}{\Lambda_{n}^{k,\mu}}
\newcommand{\Lams}{\Lambda_{n}^{k,k-1}}
\newcommand{\lamone}{\tilde\lambda_{1}}
\newcommand{\Ss}{S_{s}}
\newcommand{\Ab}[1]{A_{\nu,#1}}
\newcommand{\cb}[1]{c_{\nu,#1}}
\newcommand{\Wlat}{\W^{k,\mu}(N)}
\newcommand{\Xsp}{X_{n}^{k,\mu}}
\newcommand{\Kh}{X_{n}^{k,k-1}}
\newcommand{\II}{\mathcal I}

\title{How good are spline spaces for the Kolmogorov $n$-width of a Sobolev ball?}
\author[1]{Peng Yang\thanks{E-mail: \texttt{pyang@uestc.edu.cn}}}
\affil[1]{School of Mathematical Sciences, University of Electronic Science and Technology of China, Chengdu, Sichuan 611731, China}
\author[2]{Zhimin Zhang\thanks{E-mail: \texttt{ag7761@wayne.edu}}}
\affil[2]{Department of Mathematics, Wayne State University, Detroit, Michigan 48202, USA}
\date{}

\begin{document}

	\maketitle

	\begin{abstract}\noindent
		This paper studies the approximation power of
		polynomial splines on a uniform one-dimensional mesh.
		The quality of the best approximation of the unit ball of $H^{s}$ by an
		$n$-dimensional subspace, measured in the $L^{2}$ norm, is quantified
		by the sup-inf of that space: for each function in the Sobolev ball, take its minimal
		approximation error over the space, and then take the supremum of these
		errors over the ball. The infimum of this
		sup-inf over all $n$-dimensional spaces is the well-known Kolmogorov
		$n$-width $d_{n}$. Choosing the
		$n$-dimensional subspace to be the piecewise polynomial space
		of degree $k$ with smoothness $C^{\mu}$
		($-1\le\mu\le k-1$) gives our sup-inf $E$.
		The approximation accuracy of spline spaces is then examined by
		comparing $d_{n}$ and $E$.
		The two quantities are determined by two eigenvalue problems, through
		which we obtain a lower bound for $d_{n}$ and an upper bound for $E$.
		Under the restriction $k\ge s-1$, the upper bound for the error of the maximally smooth spline space ($\mu=k-1$) is asymptotically equal to the lower bound for $d_n$,
		showing that uniform $B$-splines are asymptotically optimal among all
		$n$-dimensional spaces. For $\mu<k-1$, we prove that $E/d_n$ is bounded below by a constant strictly larger than $1$, except when $k>s-1$ and $\mu=k-2$ hold simultaneously. This
		reveals the superiority of $B$-splines in the sense of $n$-widths.
		Numerical experiments corroborate these bounds and, in some cases
		(e.g. the exception above), show that they are sharp.
	\end{abstract}

	\section{Introduction}\label{sec:intro}

	Splines are piecewise polynomials of a prescribed degree on a mesh, glued
	at the knots with a prescribed continuity. The classical $B$-spline spaces
	are the maximally smooth choice of this construction, and they are the
	standard tool for the representation and approximation of curves and
	surfaces in computer-aided design (CAD), in the form of $B$-splines and of
	their rational counterparts, the non-uniform rational $B$-splines (NURBS)
	\cite{deBoor01,PieglTiller97,Schumaker07}.
	Spline spaces serve equally well as approximation spaces in the numerical
	treatment of differential equations \cite{LMS18}, where the smoothness of the
	space can be chosen independently of the degree. Since a CAD geometry is already described by splines, the space that
	represents it can also be used to discretize the equation posed on it.
	This observation led Hughes, Cottrell and Bazilevs to introduce
	isogeometric analysis (IGA) \cite{Cottrell09,HCB05}, in which the design
	spline space is taken directly as the space of the numerical method.
	The advantages of doing so go beyond the obvious ones: the geometry is
	represented exactly, and the remeshing step of the classical
	CAD-to-analysis workflow disappears; the smoothness of the space
	allows higher-order equations, such as the fourth-order problems of
	structural mechanics, to be discretized directly. There is, however, a further advantage that concerns the approximation itself
	\cite{BS19,EBH09}: among the spline spaces of a given polynomial degree,
	the one of maximal smoothness is the most accurate per degree of freedom.
	That property is the subject of this paper.

	The advantage of smoothness is a matter of constants rather than of rates.
	For example, on a quasi-uniform mesh in one dimension, the maximally smooth
	splines of degree $k$ and the $C^{0}$ splines of the same degree converge at
	the same asymptotic rate. Once $k\ge s-1$, the best approximation $u_{h}$ of
	a function $u\in H^{s}$ from either space satisfies
	\begin{equation}\label{eq:introrate}
		\norm{u-u_{h}}_{L^{2}}\le C\,n^{-s}\norm{u}_{H^{s}},
	\end{equation}
	where $n$ is the number of degrees of freedom; in terms of the mesh size
	$h$, this is the familiar bound of order $h^{s}$. What differs is the
	constant $C$. At the same number of degrees of freedom $n$, the numerical
	evidence consistently shows that the smoother space is the more accurate of
	the two, that is, its constant $C$ in \eqref{eq:introrate} is smaller
	\cite{EBH09,SMS20}. Error estimates with explicit constants, valid for
	spline spaces of arbitrary smoothness, have been obtained only relatively
	recently \cite{SMS20,SMS22}. For the maximally smooth space such estimates
	are given in \cite{TakacsTakacs16,SMS19}, and \cite{BeiraoVeiga11} provides
	estimates for $h$--$p$--$k$ refinement. These estimates quantify the
	advantage of smoothness observed in the simulations. However, the sharp value of the constant in \eqref{eq:introrate}
	and how it compares with the approximation constants
	of all $n$-dimensional subspaces have not been fully determined.

	To compare spaces of different smoothness at equal dimension, one needs a
	yardstick that is independent of the spaces themselves, and the natural one
	is the Kolmogorov $n$-width \cite{Kolmogorov36,Pinkus85}.
	Let $X$ be a normed linear space with norm $\norm{\cdot}_{X}$, let
	$X_{n}$ be any $n$-dimensional subspace of $X$, and let
	$\mathcal A\subset X$ be a class of functions. The sup-inf
	\begin{equation*}
		E(\mathcal A,X_{n};X):=\sup_{u\in\mathcal A}\ \inf_{u_{h}\in X_{n}}\norm{u-u_{h}}_{X}
	\end{equation*}
	measures how well $X_{n}$ approximates the worst member of $\mathcal A$ in
	the norm of $X$, and the Kolmogorov $n$-width
	\begin{equation}\label{eq:introdn}
		d_{n}(\mathcal A;X):=\inf_{\substack{X_{n}\subset X\\ \dim X_{n}=n}}E(\mathcal A,X_{n};X)
	\end{equation}
	is the least error that any $n$-dimensional subspace of $X$ can achieve.
	In this paper the general definition is specialized as follows. The
	ambient space is $X=L^{2}(0,1)$, the class is the unit ball of $H^{s}(0,1)$,
	\begin{equation*}
		\mathcal A=B^{s}(0,1):=\{u\in H^{s}(0,1):\ \norm{u}_{H^{s}(0,1)}\le1\},
	\end{equation*}
	and the spaces that are compared with $d_{n}$ are the spline spaces
	$\Xsp$ on a uniform mesh with $N$ elements. Here $k$ is the polynomial
	degree, $\mu$ is the order of continuity $C^{\mu}$ ($-1\leq \mu\leq k-1$, and $\mu=-1$ means no continuity) imposed at the knots,
	and $n=N\nu+\mu+1$, with $\nu=k-\mu\ge1$, is the dimension of the space.
	Then $E(\Bball,\Xsp;L^{2}(0,1))$ is the exact worst-case error over the unit
	ball. Equivalently, $n^{s}E(\Bball,\Xsp;L^{2}(0,1))$ is the smallest value
	that the constant $C$ in \eqref{eq:introrate} can attain. The relevant
	number is the optimality ratio
	\begin{equation}\label{eq:introLam}
		\Lam:=\frac{E(\Bball,\Xsp;L^{2}(0,1))}{d_{n}(\Bball;L^{2}(0,1))}\ \ge\ 1,
	\end{equation}
	which cannot fall below $1$ since the spline space is one of the
	competitors in \eqref{eq:introdn}, and which equals $1$ exactly when the
	spline space is itself an optimal $n$-dimensional space. The comparison of
	two $n$-dimensional spaces by the ratio of their sup-infs, and of a space
	against the width, was proposed by Babu\v{s}ka, Banerjee and Osborn
	\cite{BBO02} as a criterion for the selection of shape functions in
	generalized finite element methods. Subsequently, Evans, Bazilevs,
	Babu\v{s}ka and Hughes \cite{EBH09} applied this framework to
	one-dimensional spline spaces of variable order and continuity. For a
	suitable knot vector, the $B$-spline space is an optimal subspace for the
	$n$-width of a Sobolev-seminorm ball (a classical result of Melkman and
	Micchelli \cite{MM78}). The optimal knot locations are, however, generally
	hard to obtain.
	The numerical
	experiments in \cite{EBH09} show, perhaps not unexpectedly, that choosing a uniform knot vector leads
	to nearly optimal approximation properties of the resulting $B$-spline
	spaces.
	The following conjecture, supported by their extensive numerical
	experiments, states the asymptotic optimality of the maximally
	smooth spline space.

	\begin{conjecture}[{\cite[Conjecture 5.1]{EBH09}}]\label{conj:main}
		Let $k\ge s-1$. Then $\Lams \to 1 $ as $n\to \infty$.
	\end{conjecture}

	The restriction $k\ge s-1$ under which the conjecture is formulated is
	necessary. For $k\le s-2$, the width can be much smaller than the error of
	any piecewise polynomial space of degree $k$. Indeed, in the definition
	\eqref{eq:introdn} one may choose any $n$-dimensional space, and the width
	decays at the rate $n^{-s}$, whereas a piecewise polynomial space of degree
	$k$ on a quasi-uniform mesh cannot approximate the ball faster than
	$n^{-(k+1)}$. If, for instance, $s=k+2$, the ratio cannot tend to $1$.
	For $k\ge s-1$, no such competitor is available at the same degree, and
	Conjecture~\ref{conj:main} says that, at a
	fixed number of degrees of freedom, the maximally smooth $B$-spline space is
	asymptotically as good as the best possible space. Perhaps surprisingly, it
	also says that when the smoothness $s$ of the ball is lower than the degree
	$k$ of the space, approximating with the smoother space does not merely
	remain competitive: it is still the best choice available at the same
	number of degrees of freedom. The authors of
	\cite{EBH09} emphasize the importance of the question:
	``\emph{We feel a further analysis of the conjecture warrants merit, and
		such an analysis might prove to be instrumental to advancing spline-based
		finite element approximation theory.}'' The conjecture is settled for $s=1$, where it follows
	from the sharp one-dimensional Poincar\'e-type estimate with constant
	$h/\pi$ \cite{SMS19,SMS20} together with the explicit Neumann spectrum of
	the underlying free problem. In the variant in which the $s$-th derivative rather than the full norm is measured, optimal spline spaces have been analyzed for higher degrees and for boundary conditions by Floater and Sande \cite{FS17,FS19a,FS19b}. These optimal subspaces have been applied to
	remove spurious outliers in isogeometric eigenvalue computations
	\cite{MSpeleers22,MSpeleers23}. For general $s$, and in particular for
	$s=k+1$, which is precisely the regularity required by the standard error
	analysis for the space to attain its full order, the conjecture has
	remained open in the setting where no boundary conditions are imposed.

	Conjecture~\ref{conj:main} concerns the maximally smooth space only, and
	spaces of lower smoothness have not been analyzed in the $n$-width
	framework. That range deserves the same attention. To establish not only
	the optimality of the maximally smooth space but also its superiority over
	the spline spaces of lower smoothness, the two have to be compared within
	the same framework, and the $n$-width analysis of lower-smoothness spaces is therefore equally important.
	An existing comparison \cite{BS19} does compare spaces of different
	smoothness, but not through the width. Since the numerical evidence is
	that the smoother space is generally more accurate per degree of
	freedom, it is natural to ask whether the ratio is strictly above $1$ for every lower level of
	smoothness. For fixed integers $k\ge s-1$ and $\mu<k-1$, is $\Lam$ bounded away from $1$ as $n\to\infty$?

	We analyze both questions and obtain the following answers.

	\begingroup
	\renewcommand{\thetheorem}{A}
	\begin{theorem}\label{thm:introA}
		Let $k\ge s-1$ and $1\leq \nu\leq k+1$ with $\nu=k-\mu$. Then,
		\begin{equation*}
			\Lam\ \le\ \nu^{s}\Bigl(1+\frac{s(\mu+1)}{n}+O(n^{-2})\Bigr).
		\end{equation*}
		In particular, for the maximally smooth space ($\nu=1$),
		\begin{equation*}
			\Lams\ \le\ 1+\frac{sk}{n}+O(n^{-2}).
		\end{equation*}
		Consequently, $ \Lams \to 1 $ as $n \to \infty$.
	\end{theorem}
	\endgroup

	\begingroup
	\renewcommand{\thetheorem}{B}
	\begin{theorem}\label{thm:introB}
		Let $k\ge s-1$ and $1\leq \nu\leq k+1$ with $\nu=k-\mu$. Then,
		\begin{equation*}
			\Lam
			\ \ge\ \nu^{s}\,\Bigl(\frac{\Ab{k+1-s}}{\Ab{k+1}}\Bigr)^{1/2}\,\bigl(1-o(1)\bigr),
		\end{equation*}
		where
		\begin{equation*}
			\Ab{r}:=4^{-\nu}\sum_{m=0}^{\nu}\binom{\nu}{m}^{2}\abs{\nu-2m}^{2r},
		\end{equation*}
		and $0^0=1$.
		Consequently, for a sufficiently large $n$,
		\begin{enumerate}[label=(\roman*),leftmargin=2.4em]
				\item If $k=s-1$, then $\Lam>1$ for every $\nu\ge2$;
			\item If $k>s-1$, then $\Lam>1$ for every $\nu>2$, while for $\nu=2$, the constant factor in the lower bound equals $1$.
		\end{enumerate}
	\end{theorem}
	\endgroup

	The two theorems are proved by analyzing the numerator and denominator of the ratio \eqref{eq:introLam} separately; each is governed by an eigenvalue problem.
	For the $n$-width, a
	crude lower bound suffices. The main task is the sup-inf, which we have found
	to be the square root of the largest quotient of a squared $\Hm$-norm by a
	squared $L^{2}$-norm, taken over the $L^{2}$-orthogonal complement of the
	spline space.
	The complement has an explicit description. The $(k+1)$-st antiderivative
	of a function in the orthogonal complement is determined uniquely: it vanishes at every
	interior knot together with several of its derivatives, their number
	depending on $\mu$, and at the endpoints together with all its derivatives
	up to order $k$. The vanishing forces these antiderivatives to oscillate,
	so that the mesh frequency appears in their derivatives. Everything is
	thus reduced to the $L^{2}$ and $\Hm$ norms of the antiderivatives.
	An upper bound for the $\Hm$-norm gives a crude upper bound for the sup-inf.
	Although not sharp for all spline spaces, it already shows that the maximally
	smooth space is asymptotically optimal, in the sense that its error asymptotically matches the Kolmogorov $n$-width
	as $n\to\infty$. This yields the first theorem.
	Lower bounds are obtained by testing the $\Hm$-norm against two different
	functions, a single oscillating function and a cosine sequence. The first
	choice is elementary. As a lower bound for the ratio, it yields information
	on the comparison of a lower-smoothness space with the maximally smooth
	spline space only when $\nu\ge3$. The second choice requires a finer analysis,
	which yields the second theorem above. The numerical section corroborates
	the theory and shows that the lower bound in Theorem~\ref{thm:introB} is accurate. When the two bounds obtained from the two test functions
	nearly coincide, the larger one, the cosine-sequence bound, is approached from
	below by the computed values; these computed values serve as our reference for the true value, and their Galerkin error is negligible. This
	indicates that the bound is sharp. When they are far apart, the
	cosine-sequence bound is by far the more accurate of the two and remains
	close to the computed values.

	The paper is organized as follows. Section~\ref{sec:setting} gives a lower
	bound for the $n$-width and shows that the sup-inf depends only on the
	$L^{2}$-orthogonal complement of the spline space, for which it gives an explicit description.
	Section~\ref{sec:floor} proves a crude upper bound for the
	sup-inf and deduces Theorem~\ref{thm:introA}. Section~\ref{sec:nu} tests
	the $\Hm$-norm against two functions in $\Hsp$, obtains a different lower
	bound for the sup-inf from each, and proves Theorem~\ref{thm:introB}.
	Section~\ref{sec:numerics} reports the numerical experiments and shows
	that the lower bound is sharp when the two bounds coincide.

	Throughout the paper, $O(a_n)$ with $a_n>0$ denotes a quantity bounded in absolute value by a constant multiple of $a_n$. In
	an identity, $f=g+O(a_n)$ means $\abs{f-g}\le c\,a_n$ for a constant $c$
	independent of the parameters. In an inequality, $O(a_n)$ denotes a positive
	quantity: $f\leq g+O(a_n)$ means $f-g\le c\,a_n$, and $f\geq g-O(a_n)$
	means $g-f\leq c\,a_n$. Consequently, $f=g+O(a_n)$ denotes
	$f\leq g+O(a_n)$ and $f\geq g-O(a_n)$ simultaneously. Moreover,
	$b_n=o(a_n)$ means that $\abs{b_n/a_n}\to0$ as $n\to\infty$. When
	$o(a_n)$ appears in an inequality, it denotes a positive quantity as well,
	with the same convention as $O(a_n)$ above.

	\section{The \texorpdfstring{$n$}{n}-width and the sup-inf}\label{sec:setting}

	This section fixes the notation and records the two characterizations on
	which everything below rests. The first is the value of the Kolmogorov
	$n$-width and an optimal subspace for it, which come from the eigenpairs of
	the free problem. The second expresses the sup-inf of a spline space through
	the $\Hm$-norm on the orthogonal complement of the space.

	\subsection{The Kolmogorov \texorpdfstring{$n$}{n}-width}\label{sec:notations}

	Let $\Omega=(0,1)$, and let $\Hsp=\Hsp(\Omega)$, $s\ge1$, carry the
	full inner product
	\begin{equation*}
		(u,v)_{\Hsp}=\sum_{a=0}^{s}\int_{0}^{1}u^{(a)}v^{(a)}\,dx,
		\qquad
		\norm{u}^{2}_{\Hsp}=\sum_{a=0}^{s}\norm{u^{(a)}}^{2}_{L^{2}}.
	\end{equation*}
	All functions are real-valued and all
	spaces are real Hilbert spaces, and
	\begin{equation*}
		\Bball:=\{u\in\Hsp:\ \norm{u}_{\Hsp}\le1\}
	\end{equation*}
	is the unit ball of $\Hsp$. The worst-case error of a space $X\subset L^{2}$ is
	\begin{equation}\label{eq:sup-infs}
		E(\Bball,X;L^{2}):=\sup_{u\in\Bball}\ \inf_{v\in X}\norm{u-v}_{L^{2}},
	\end{equation}
	and the best any $n$-dimensional space can do is
	\begin{equation}\label{eq:widths}
		\dn:=\inf_{\dim X_{n}=n}E(\Bball,X_{n};L^{2}),
	\end{equation}
	which is the Kolmogorov $n$-width of the ball \cite{Kolmogorov36}.

	The width is an infimum over all $n$-dimensional
	subspaces of $L^{2}$, and both its value and an optimal subspace are
	spectral. The spectrum in question is that of the free problem
	\begin{equation}\label{eq:free}
		(u,v)_{L^{2}}=\lambda\,(u,v)_{\Hsp}\qquad\text{for all }v\in\Hsp,
	\end{equation}
	in which no boundary conditions are imposed. The associated operator is self-adjoint with compact
	resolvent, because the embedding $\Hsp\hookrightarrow L^{2}$ is compact, so
	it has eigenvalues $\lambda_{1}\ge\lambda_{2}\ge\cdots>0$
	and eigenfunctions $\varphi_{1},\varphi_{2},\dots$ forming an orthonormal
	basis of $L^{2}$. We normalize them by
	\begin{equation}\label{eq:normphi}
		(\varphi_{j},\varphi_{k})_{L^{2}}=\delta_{jk},
		\qquad\text{whence}\qquad
		(\varphi_{j},\varphi_{k})_{\Hsp}=\lambda_{j}^{-1}\delta_{jk},
	\end{equation}
	the second relation being \eqref{eq:free} with $u=\varphi_{j}$ and $v=\varphi_{k}$.
	For an eigenfunction $u=\varphi_{j}$, \eqref{eq:free} with $v=u$ gives
	\[
	\rho_{j}:=\lambda_{j}^{-1}
	=\frac{\norm{\varphi_{j}}^{2}_{\Hsp}}{\norm{\varphi_{j}}^{2}_{L^{2}}} .
	\]
	Since $\norm{u}_{\Hsp}\ge\norm{u}_{L^{2}}$ for every $u\in\Hsp$, this gives
	$\rho_{j}\ge1$ for every $j$, and $\rho_{1}=1$ since the constant function is a
	free eigenfunction. The reciprocal eigenvalues therefore increase from
	$\rho_{1}=1$, that is,
	\begin{equation}\label{eq:rhoorder}
		1=\rho_{1}\le\rho_{2}\le\rho_{3}\le\cdots.
	\end{equation}
	From this expression, we immediately have the following lemma characterizing the Kolmogorov $n$-width of the ball.
	\begin{lemma}[{\cite[Theorem~4.5]{EBH09}}]\label{lem:dneig}
		For every $n\ge1$,
			$\dn^{2}=\lambda_{n+1}=\rho_{n+1}^{-1}$.
	\end{lemma}

	For the width, we shall use the following rough bound, together with the
	min--max form of the reciprocal eigenvalues.

	\begin{lemma}\label{lem:minmax}
		For every $n\ge1$, the reciprocal eigenvalues of the free problem \eqref{eq:free} satisfy
		\begin{equation}\label{eq:minmax}
			\rho_{n+1}=\min_{\dim U=n+1}\ \max_{0\ne u\in U}
			\frac{\norm{u}^{2}_{\Hsp}}{\norm{u}^{2}_{L^{2}}},
		\end{equation}
		and
		\begin{equation}\label{eq:bound-rho}
			\rho_{n+1} \leq \sum_{a=0}^{s}(n\pi)^{2a}.
		\end{equation}
		Consequently, we have a lower bound for the Kolmogorov $n$-width
		\begin{equation}\label{eq:crude}
			\dn=\rho_{n+1}^{-1/2}\ge(n\pi)^{-s}\bigl(1-O(n^{-2})\bigr).
		\end{equation}
	\end{lemma}

	\begin{proof}
		For $u\in H^s$, write
		$u=\sum_{j\ge1}c_j\varphi_j$ with $c_j=(u,\varphi_j)_{L^2}$. The orthogonality
		relations of \eqref{eq:normphi} give
		\begin{equation*}
			\|u\|^2_{L^2}=\sum_jc_j^2, \qquad \|u\|^2_{H^s}=\sum_j\rho_jc_j^2.
		\end{equation*}
		Define
		\begin{equation*}
			R(u):=\frac{\|u\|^2_{H^s}}{\|u\|^2_{L^{2}}}
			=\frac{\sum_j\rho_jc_j^2}{\sum_jc_j^2}.
		\end{equation*}
		Choosing $U=\mathrm{span}\{\varphi_1,\dots,\varphi_{n+1}\}$ gives
		$\max_{0\ne u\in U}R(u)=\rho_{n+1}$ (since $\{\rho_{j}\}_{j\geq 1}$ increases), so
		\begin{equation*}
			\min_{\dim U=n+1}\ \max_{0\ne u\in U}
			\frac{\norm{u}^{2}_{\Hsp}}{\norm{u}^{2}_{L^{2}}}\leq \rho_{n+1}.
		\end{equation*}
		Conversely, let $U$ be any $(n+1)$-dimensional subspace and
		put $V:=\mathrm{span}\{\varphi_1,\dots,\varphi_n\}$. Pick $0\ne u\in U\cap V^{\perp}$ ($\perp$ denoting the $L^2$-orthogonal complement),
		i.e.\ $c_j=0$ for $j\le n$. Then $R(u)=\sum_{j\geq n+1}\rho_jc_j^2/\sum_{j\geq n+1}c_j^2\ge\rho_{n+1}$, so
		$\max_{0\ne u\in U}R(u)\ge \max_{0\ne u\in U\cap V^{\perp}}R(u) \geq \rho_{n+1}$. Since $U$ was an arbitrary
		$(n+1)$-dimensional space, passing to the infimum gives
		\begin{equation*}
			\min_{\dim U=n+1}\ \max_{0\ne u\in U}
			\frac{\norm{u}^{2}_{\Hsp}}{\norm{u}^{2}_{L^{2}}}\geq \rho_{n+1}.
		\end{equation*}
		This proves \eqref{eq:minmax}.

		To show \eqref{eq:bound-rho}, we consider the $(n+1)$-dimensional space
		\begin{equation*}
			U= \mathrm{span}\{1,\cos(\pi x),\dots,\cos(n\pi x)\}.
		\end{equation*}
		The basis functions are mutually orthogonal on
		$(0,1)$, and so are their derivatives: for $a\ge1$ the functions
		$\frac{d^a}{dx^a}\cos(k\pi x)$, $1\le k\le n$, are sines for odd $a$ and cosines for
		even $a$, all with the frequencies $k\pi$, hence mutually orthogonal, and the
		constant function contributes nothing to the derivatives. Alternating between sines and cosines changes
		no norm, because $\frac{d^a}{dx^a}\cos(k\pi x)=(k\pi)^a\cos(k\pi x+a\pi/2)$
		and, for integer $k\ge1$ and every phase $\theta$,
		\[
		\|\cos(k\pi\cdot+\theta)\|^2_{L^{2}}=\int_0^1\cos^2(k\pi x+\theta)\,dx
		=\int_0^1\frac{1+\cos(2k\pi x+2\theta)}{2}\,dx=\frac12.
		\]
		Hence
		\[
		\Bigl\|\frac{d^a}{dx^a}\cos(k\pi x)\Bigr\|^2_{L^{2}}
		=(k\pi)^{2a}\|\cos(k\pi\cdot)\|^2 _{L^{2}}\qquad(a\ge0,\ 1\le k\le n),
		\]
		where the factor $\|\cos(k\pi\cdot)\|^2_{L^{2}}=1/2$ is the same constant for sines ($a$ odd) and cosines ($a$ even). For any $u\in U$, write it as
		$u=c_0+\sum_{1\le k\le n}c_k\cos(k\pi x)$, which satisfies, for every $a\ge1$,
		\[
		\|u^{(a)}\|^2_{L^{2}}=\sum_{1\le k\le n}(k\pi)^{2a}c_k^2\|\cos(k\pi\cdot)\|^2_{L^{2}}
		\le(n\pi)^{2a}\|u\|^2_{L^{2}},
		\]
		while for $a=0$ the same inequality holds, with $(n\pi)^0=1$. Therefore,
		\begin{equation*}
			\|u\|^2_{H^s}=\sum_{0\le a\le s}\|u^{(a)}\|^2_{L^{2}}\le \sum_{a=0}^{s}(n\pi)^{2a}\|u\|^2_{L^{2}}.
		\end{equation*}
		This completes the proof.
	\end{proof}

	The crude bound \eqref{eq:crude} is not sharp, and the constant that replaces
	it is set by the spectrum of the free problem, which we now describe. It is
	convenient to write \eqref{eq:free} in strong form. Integrating
	$\int_{0}^{1}u^{(a)}v^{(a)}$ by parts $a$ times in each summand gives, for
	$u\in H^{2s}(0,1)$ and $v\in\Hsp$,
	\begin{equation*}
		(u,v)_{\Hsp}=\int_{0}^{1}(\mathcal L u)\,v\,dx
		+\sum_{m=0}^{s-1}\big((\mathcal B_mu)v^{(m)}\big)(x)\big|_{x\in \{0,1\}}\,,
	\end{equation*}
	with operators
	\begin{equation*}
		\mathcal L:=\sum_{a=0}^{s}(-1)^{a}\frac{d^{2a}}{dx^{2a}},\qquad \mathcal B_m:= \sum_{a=m+1}^{s}(-1)^{a-1-m}\frac{d^{2a-1-m}}{dx^{2a-1-m}}.
	\end{equation*}
	The values
	$v^{(0)},\dots,v^{(s-1)}$ may be prescribed arbitrarily at the two endpoints,
	so the free problem \eqref{eq:free} is equivalent to the eigenvalue problem
	\begin{equation}\label{eq:strong}
		\mathcal L u=\rho u\ \ \text{on }(0,1),
		\qquad
		(\mathcal B_mu)(x)\big|_{x\in\{0,1\}}=0,\quad m=0,\dots,s-1.
	\end{equation}
	We first give the sharp estimates of the eigenvalues in the two lowest
	cases, $s=1$ and $s=2$; for general $s$, we then state the corresponding sharp
	form as a conjecture.
	\begin{lemma}\label{lem:spectrum}
		For $s=1,2$, the eigenvalues of \eqref{eq:strong} satisfy
		\begin{equation}\label{eq:sharp-estimate}
		\rho_{n+1}^{1/(2s)}=\pi\bigl(n-\frac{s-1}{2}\bigr)+O(n^{-1}).
		\end{equation}
	\end{lemma}

	\begin{proof}
		By \eqref{eq:strong} with $s=1$, we have
		\begin{equation*}
			-u''+u=\rho u\ \ \text{on }(0,1), \qquad u'(x)\big|_{x\in\{0,1\}}=0.
		\end{equation*}
		Its eigenfunctions are the cosines $\cos((j-1)\pi x)$ and
		$\rho_{j}=1+((j-1)\pi)^{2}$. The desired result \eqref{eq:sharp-estimate} for $s=1$ is then clear.

		For \eqref{eq:strong} with $s=2$, we have
		\begin{equation}\label{eq:s2-eigenvalue}
			u^{(4)}-u''+(1-\rho)u=0\ \ \text{on }(0,1), \qquad u''(x)\big|_{x\in\{0,1\}}=0,\,\,(u'''-u')(x)\big|_{x\in\{0,1\}}=0.
		\end{equation}
		Consider the characteristic equation
		\begin{equation*}
			r^{4}-r^{2}+1=\rho,
		\end{equation*}
		which gives $r^2 = (1\pm \sqrt{4\rho-3})/2$. Now parametrize the eigenvalue by a single number $p$ instead of by $\rho$: let
		\begin{equation}\label{eq:rho-p}
			\rho = 1+p^2+p^4, \quad p\geq 0.
		\end{equation}
		This loses nothing: every eigenvalue is $\ge1$ by \eqref{eq:rhoorder}, and every
		$\rho\ge1$ is $1+p^{2}+p^{4}$ for exactly one $p\ge0$. Moreover,
		$4\rho-3=4(1+p^{2}+p^{4})-3=(1+2p^{2})^{2}$, so
		\begin{equation*}
			r^2 = 1+p^2=:q^2 \,\,(q>0)\,\, \text{or } r^2 = -p^2, \quad \text{i.e. } \quad r = \pm q\,\, \text{or } r = \pm p\, i.
		\end{equation*}
		The general solution of \eqref{eq:s2-eigenvalue} is then
		\begin{equation*}
			u(x) = c_1 \,\sin(px) + c_2 \,\cos(px) + c_3 \,\sinh(qx) + c_4 \,\cosh(qx).
		\end{equation*}
		Note that both the equation and the conditions in \eqref{eq:s2-eigenvalue}
		are invariant under $x\mapsto1-x$: if $u$ solves the problem then so does
		$Ru(x):=u(1-x)$. Hence $u+Ru$, which is even about $x=\frac12$, and $u-Ru$,
		which is odd, also solve the problem, so every solution is the sum of an even
		and an odd solution, and it suffices to treat the two classes separately.
		Since $\cos(p(x-\frac12))$ and $\cosh(q(x-\frac12))$ are even about
		$x=\frac12$, while $\sin(p(x-\frac12))$ and $\sinh(q(x-\frac12))$ are odd,
		these classes are
		\[
		u(x)=A\cos\Bigl(p\bigl(x-\tfrac12\bigr)\Bigr)+B\cosh\Bigl(q\bigl(x-\tfrac12\bigr)\Bigr),
		\qquad
		u(x)=C\sin\Bigl(p\bigl(x-\tfrac12\bigr)\Bigr)+D\sinh\Bigl(q\bigl(x-\tfrac12\bigr)\Bigr).
		\]
		For $u$ in the even class, the boundary conditions in
		\eqref{eq:s2-eigenvalue} read
		\[
		-Ap^{2}\cos\frac p2+Bq^{2}\cosh\frac q2=0,\qquad
		Aq\sin\frac p2+Bp\sinh\frac q2=0,
		\]
		where the second equation uses the relation $1+p^{2}=q^{2}$. The eigenfunctions for
		the parameter $p$ are the nonzero solutions of $M(p)\binom AB=0$, with
		$M(p)$ the matrix of these two equations, and the map $(A,B)\mapsto u$ is injective. A repeated eigenvalue in the class would mean
		$\dim\ker M(p)\ge2$, that is, $M(p)=\boldsymbol{0}$; but the
		entry $q^{2}\cosh\frac q2$ never vanishes, because $q=\sqrt{1+p^{2}}\ge1$ and
		$\cosh>0$. So an eigenvalue of the even class is simple. The determinant of
		$M(p)$ vanishes exactly when
		\begin{equation}\label{eq:s2det1}
			\tan\frac p2=-\frac{p^{3}}{q^{3}}\tanh\frac q2 .
		\end{equation}
		The odd class is treated in the same way. The
		boundary conditions read
		\[
		-Cp^{2}\sin\frac p2+Dq^{2}\sinh\frac q2=0,\qquad
		-Cq\cos\frac p2+Dp\cosh\frac q2=0.
		\]
		The same argument
		as above shows that an eigenvalue of the odd class is simple, and a nonzero solution exists if and only if
		\begin{equation}\label{eq:s2det2}
			\tan\frac p2=\frac{q^{3}}{p^{3}}\tanh\frac q2.
		\end{equation}
		Clearly, the two equations \eqref{eq:s2det1} and \eqref{eq:s2det2} do not share a root in $p$.

		The eigenvalues are therefore the parameters $p\ge0$ that solve
		\eqref{eq:s2det1} or \eqref{eq:s2det2}, and both equations are analyzed with
		the two functions
		\[
		h(p):=-\frac{p^{3}}{q^{3}}\tanh\frac q2 \qquad(p\geq0),
		\qquad\quad
		f(p):=\frac{q^{3}}{p^{3}}\tanh\frac q2
		\qquad(p>0),
		\]
		which are monotone: differentiating, one finds that both are decreasing. Moreover, $h(p)\to0$, $f(p)\to+\infty$ as $p\to0^{+}$,
		and $h(p)\to -1$, $f(p)\to 1$ as $p\to\infty$.
		Hence $h$ decreases from $0$ to $-1$, and $f$ decreases from $+\infty$ to
		$1$.
		Now note that $\tan p/2$ increases from $-\infty$ to $0$ on each
		$\bigl(\pi+2\pi m,2\pi+2\pi m\bigr)$ and from $1$ to $+\infty$ on each interval
		$\bigl(\frac\pi2+2\pi m,\pi+2\pi m\bigr)$, $m\ge0$. Consequently,
		\begin{itemize}
			\item on each $\bigl(\pi+2\pi m,2\pi+2\pi m\bigr)$, the even equation
			\eqref{eq:s2det1} has exactly one solution, at which
			$\tan\frac p2=-1+O(m^{-2})$, that is, $p=\frac{3\pi}2+2\pi m+O(m^{-2})$;
			\item on each $\bigl(\frac\pi2+2\pi m,\pi+2\pi m\bigr)$, the odd equation
			\eqref{eq:s2det2} has exactly one solution, at which
			$\tan\frac p2=1+O(m^{-2})$, that is, $p=\frac\pi2+2\pi m+O(m^{-2})$;
			\item on each $\bigl(2\pi m,\frac\pi2+2\pi m\bigr)$, $m\ge0$, neither equation has a solution; and $p=0$ solves
			\eqref{eq:s2det1} and yields the constant eigenfunction, of eigenvalue
			$\rho=1$.
		\end{itemize}
		Listing the eigenvalues in increasing order gives $p_{1}=0$ and
		$p_{j}=\pi\bigl(j-\frac32\bigr)+O(j^{-2})$ for $j\ge2$, so that, by \eqref{eq:rho-p},
		\[
		\rho_{j}^{1/4}=p_{j}\bigl(1+p_{j}^{-2}+p_{j}^{-4}\bigr)^{1/4}
		=p_{j}+\frac{1}{4p_{j}}+O(p_{j}^{-3})
		=\pi\bigl(j-\tfrac32\bigr)+O(j^{-1}).
		\]
		This is the case $s=2$ of \eqref{eq:sharp-estimate}, which completes the proof.
	\end{proof}

	\begin{conjecture}\label{conj:shift}
		For every $s\ge1$, the eigenvalues of \eqref{eq:strong} satisfy
		\begin{equation}\label{eq:shiftC}
			\rho_{n+1}^{1/(2s)}=\pi\bigl(n-\frac{s-1}{2}\bigr)+O(n^{-1}).
		\end{equation}
	\end{conjecture}

	The same asymptotic is conjectured by Bressan, Floater and Sande \cite{BFS21} for
	the seminorm problem, in which the lower-order terms of $\mathcal L$ are dropped
	and the operator is the pure $(-1)^{s}D^{2s}$.
	Conjecture~\ref{conj:shift} is illustrated numerically in
	Section~\ref{sec:numerics} (Table~\ref{tab:1}).

	\subsection{The sup-inf of a spline space}\label{sec:twoquantities}

	We now take the space $X$ in the sup-inf \eqref{eq:sup-infs} to be an
	$n$-dimensional space of piecewise polynomials on a uniform mesh, with possibly high smoothness, i.e. a spline space. Then,
	the sup-inf \eqref{eq:sup-infs} is the largest error, measured in $L^{2}$, that the
	best approximation from the spline space can leave on the ball $B^s$.

	For a degree $k\ge0$ and an integer $N\ge k+1$, the mesh is uniform with its knots $t_{p}=ph$,
	$h=1/N$, $p=0,\dots,N$, so that the elements are the $N$ intervals
	$(t_{p},t_{p+1})$, all of the same length $h$. For a smoothness $\mu\in\{-1,0,\dots,k-1\}$, let
	\begin{equation}\label{eq:X}
		\Xsp:=\bigl\{v:\ v|_{(t_{p},t_{p+1})}\in\Pk\ \ (0\le p\le N-1),
		\ \ v\in C^{\mu}(0,1)\bigr\}
	\end{equation}
	be the space of functions that agree on each element with a polynomial
	of degree $\le k$ and whose derivatives of order $\le\mu$ are continuous
	across the knots ($\mu=-1$ means that no continuity is required); here $\Pk$
	denotes the polynomials of degree $\le k$. The space has dimension
	\begin{equation}\label{eq:dimX}
		n=\dim\Xsp=N(k+1)-(N-1)(\mu+1)=N\nu+\mu+1,
	\end{equation}
	with
	\begin{equation*}
		\nu:=k-\mu\ge1.
	\end{equation*}
	Two values of the smoothness are extreme. The discontinuous space corresponds to $\mu=-1$; it has $\nu=k+1$ degrees of freedom per element and dimension $n=N(k+1)$. The other extreme is the maximally smooth space
	$\mu=k-1$, of continuity $C^{k-1}$, and dimension $n=N+k$; it is the space appearing in
	Conjecture~\ref{conj:main}.

	Let $\Pi_{n}$ be the $L^{2}$-orthogonal projection onto $\Xsp$. The
	sup-inf \eqref{eq:sup-infs} with $X= \Xsp$ is characterized by the eigenvalue problem
	\begin{equation}\label{eq:supeig}
		\bigl((I-\Pi_{n})u,\,(I-\Pi_{n})v\bigr)_{L^{2}}
		=\tilde\lambda\,(u,v)_{\Hsp}
		\qquad\text{for all }v\in\Hsp.
	\end{equation}
	Its eigenfunctions are the nonzero $u\in\Hsp$ for which \eqref{eq:supeig}
	holds for some $\tilde\lambda\in\mathbb{R}$, and the problem is well-posed,
	with at most $n$ zero eigenvalues. If $u\in\Xsp\cap\Hsp$, then
	$(I-\Pi_{n})u=0$ and \eqref{eq:supeig} holds with $\tilde\lambda=0$.
	Conversely, $\tilde\lambda=0$ forces $(I-\Pi_{n})u=0$, that is,
	$u\in\Xsp\cap\Hsp$. The nonzero
	eigenvalues form a decreasing sequence
	$\tilde\lambda_{1}\ge\tilde\lambda_{2}\ge\cdots>0$, and
	\begin{equation}\label{eq:lambda1}
		\tilde\lambda_{1}=\sup_{0\ne u\in\Hsp}
		\frac{\norm{(I-\Pi_{n})u}^{2}_{L^{2}}}{\norm{u}^{2}_{\Hsp}}
	\end{equation}
	is their largest element. The eigenvalue $\tilde\lambda_{1}$ is the
	relative remainder left by the projection, in the worst case over
	$u\in\Hsp$. The sup-inf of the space is its square root,
	\begin{equation}\label{eq:Etilde}
		E(\Bball,\Xsp;L^{2})=\tilde\lambda_{1}^{1/2},
	\end{equation}
	because \eqref{eq:sup-infs} may be evaluated at the projection: for every
	$u$ the two summands of $u=\Pi_{n}u+(I-\Pi_{n})u$ are $L^{2}$-orthogonal, so
	that $\mathrm{dist}_{L^{2}}(u,\Xsp)=\norm{(I-\Pi_{n})u}_{L^{2}}$.

	By the Riesz representation
	theorem, for every $g\in L^{2}$ there is a unique $z_{g}\in\Hsp$ with
	\begin{equation}\label{eq:Bdef}
		(z_{g},v)_{\Hsp}=(g,v)_{L^{2}}\qquad\text{for all }v\in\Hsp ,
	\end{equation}
	and we define $Rg:=z_{g}$. The resulting operator $R:L^{2}\to\Hsp$ is
	linear and bounded. The lemma below gives another form of the largest
	eigenvalue $\lamone$.

	\begin{lemma}\label{lem:collapse}
		Write $P:=I-\Pi_{n}$, and let $R$ be the Riesz map \eqref{eq:Bdef}. Then $R^{1/2}$ is an isometry of $L^{2}$ onto $\Hsp$, and
		\begin{equation}\label{eq:collapse}
			\lamone=\norm{PR^{1/2}}^{2}=\norm{R^{1/2}P}^{2}
			=\sup_{0\ne f\in(\Xsp)^{\perp}}\frac{\norm{f}^{2}_{\Hm}}{\norm{f}^{2}_{L^{2}}},
		\end{equation}
		where $(\Xsp)^{\perp}$ is the $L^{2}$-orthogonal complement of $\Xsp$, and $\norm{\cdot}_{\Hm}$ is the dual norm
		\begin{equation}\label{eq:dualmnorm}
			\norm{f}_{\Hm}:=\sup_{0\ne\delta\in\Hsp}
			\frac{(\delta,f)_{L^{2}}}{\norm{\delta}_{\Hsp}}.
		\end{equation}
	\end{lemma}

	\begin{proof}
		One checks from \eqref{eq:Bdef} that $R$ is self-adjoint for both inner
		products and positive: $(Rg,h)_{L^{2}}=(g,Rh)_{L^{2}}$ for $g,h\in L^{2}$,
		$(Rf,h)_{\Hsp}=(f,Rh)_{\Hsp}$ for $f,h\in\Hsp$, and
		$(Rg,g)_{L^{2}}=\norm{Rg}^{2}_{\Hsp}\ge0$.
		As a bounded symmetric positive operator on $L^{2}$, $R$ has a bounded
		symmetric positive square root $R^{1/2}$. Moreover, $R^{1/2}$ is symmetric
		for the $\Hsp$-inner product as well, since $R$ is. Hence, for $u\in\Hsp$,
		\[
		\norm{R^{1/2}u}^{2}_{\Hsp}=(Ru,u)_{\Hsp}=\norm{u}^{2}_{L^{2}}.
		\]
		Since $\Hsp$ is dense
		in $L^{2}$, while $R^{1/2}$ is bounded on $L^{2}$ and $\Hsp$ is complete,
		the identity extends to all $u\in L^{2}$. Thus $R^{1/2}$ is an
		isometry of $L^{2}$ onto $\Hsp$. Its range is closed, being the isometric
		image of the complete space $L^{2}$, and dense, since a $w\in\Hsp$
		orthogonal to $R^{1/2}(L^{2})$ satisfies
		$0=(R^{1/2}R^{1/2}w,w)_{\Hsp}=(Rw,w)_{\Hsp}=\norm{w}^{2}_{L^{2}}$.

		Consequently, for $g\in L^{2}$ the dual norm is computed by substituting
		$\delta=R^{1/2}\eta$ in the supremum \eqref{eq:dualmnorm}. As $\eta$ runs
		through $L^{2}$, the element $\delta$ runs through all of $\Hsp$, with
		$\norm{\delta}_{\Hsp}=\norm{\eta}_{L^{2}}$ by the isometry. Then the
		supremum is the dual characterization of the $L^{2}$ norm
		\begin{equation}\label{eq:Hm}
			\norm{g}_{\Hm}
			=\sup_{0\ne\eta\in L^{2}}\frac{(R^{1/2}\eta,g)_{L^{2}}}{\norm{\eta}_{L^{2}}}=\sup_{0\ne\eta\in L^{2}}\frac{(\eta,R^{1/2}g)_{L^{2}}}{\norm{\eta}_{L^{2}}}
			=\norm{R^{1/2}g}_{L^{2}}.
		\end{equation}
		It remains to insert these identities into the characterization of
		$\lamone$. Replacing $u$ by $R^{1/2}v$ in \eqref{eq:lambda1} gives
		\[
		\lamone=\sup_{0\ne v\in L^{2}}\frac{\norm{PR^{1/2}v}^{2}_{L^{2}}}{\norm{R^{1/2}v}_{H^s}^{2}}=\sup_{0\ne v\in L^{2}}\frac{\norm{PR^{1/2}v}^{2}_{L^{2}}}{\norm{v}^{2}_{L^{2}}}
		=\norm{PR^{1/2}}^{2}.
		\]
		Since $(PR^{1/2})^{*}=R^{1/2}P$, the two operator norms in \eqref{eq:collapse} agree. By the
		definition of the operator norm,
		\[
		\|R^{1/2}P\|^2
		=\sup_{0\ne g\in L^2}\frac{\|R^{1/2}Pg\|^2_{L^2}}{\|g\|^2_{L^2}} ,
		\]
		and the numerator depends on $g$ only through $f:=Pg$. Since
		$\|g\|^2_{L^2}=\|f\|^2_{L^2}+\|(I-P)g\|^2_{L^2}$, replacing the denominator
		$\|g\|^2_{L^2}$ by the smaller $\|f\|^2_{L^2}$ can only increase the quotient:
		for every $g$ with $f=Pg\ne0$,
		\[
		\frac{\|R^{1/2}Pg\|^2_{L^2}}{\|g\|^2_{L^2}}
		\le\frac{\|R^{1/2}f\|^2_{L^2}}{\|f\|^2_{L^2}} ,
		\]
		with equality whenever $g$ itself lies in $P(L^2)$ (then $f=g$). Hence the supremum over all
		$g\ne0$ equals the supremum over the nonzero elements of the range $P(L^2)$,
		a nonzero space because $X_n$ is finite dimensional, and by \eqref{eq:Hm} this
		is
		\[
		\|R^{1/2}P\|^2
		=\sup_{0\ne f\in P(L^2)}\frac{\|R^{1/2}f\|^2_{L^2}}{\|f\|^2_{L^2}}
		=\sup_{0\ne f\in P(L^2)}\frac{\|f\|^2_{H^{-s}}}{\|f\|^2_{L^2}} .
		\]
		Finally, since $P(L^2)=(X_n^{k,\mu})^{\perp}$, the proof of \eqref{eq:collapse} is complete.
	\end{proof}

	By Lemma~\ref{lem:collapse}, $\lamone$ is the largest quotient of a squared
	$\Hm$-norm by a squared $L^{2}$-norm, over the $L^{2}$-orthogonal
	complement of the spline space. The distance to $\Xsp$ has disappeared
	from the original definition \eqref{eq:lambda1}, and what remains is the problem of describing
	that complement. At the end of this subsection, we give an
	explicit characterization of the complement $(\Xsp)^{\perp}$; it is the basis on which the
	analysis of the following sections rests.

	\begin{definition}\label{def:J}
		For $f\in L^{2}(0,1)$, let
		\begin{equation*}
			F(x) = (\J^{k+1}f)(x):=\frac{1}{k!}\int_{0}^{x}(x-r)^{k}f(r)\,dr
		\end{equation*}
		be the $(k+1)$-fold antiderivative. By construction, $F^{(k+1)}=f$ and $F^{(l)}(0)=0\ \ (l=0,\dots,k)$.
	\end{definition}

	\begin{definition}\label{def:W}
		With $\nu:=k-\mu\geq 1$,
		\begin{equation*}
			\begin{aligned}
				\Wlat:=\Bigl\{F\in H^{k+1}(0,1):\;&
				F^{(l)}(t_{p})=0\ \ (l=0,\dots,\nu-1,\ p=1,\dots,N-1),\\
				&F^{(l)}(0)=F^{(l)}(1)=0\ \ (l=0,\dots,k)\Bigr\}.
			\end{aligned}
		\end{equation*}
		The first line says that $F$ vanishes to order $\nu$ at every interior
		knot, the second that $F$ vanishes to order $k+1$ at each endpoint; for
		the maximally smooth space $\nu=1$, the first line is the single condition
		$F(t_{p})=0$.
	\end{definition}

	\begin{theorem}\label{thm:structure}
		Let $X_n^{k,\mu}$ be the spline space \eqref{eq:X}.
		For every $f\in L^{2}(0,1)$,
		\begin{equation*}
			f\in(\Xsp)^{\perp}
			\qquad\Longleftrightarrow\qquad
			\J^{k+1}f\in\Wlat .
		\end{equation*}
		Consequently, $\J^{k+1}$ is a bijection from $(\Xsp)^{\perp}$ onto $\Wlat$,
		with inverse $F\mapsto F^{(k+1)}$, and
		\begin{equation*}
			(\Xsp)^{\perp}=\{F^{(k+1)}:F\in\Wlat\},
		\end{equation*}
		where the derivative is understood in the weak sense.
	\end{theorem}

	\begin{proof}
		Note that $f\in(\Xsp)^{\perp}$ means
		$(f,v)_{L^2}=0$ for every $v\in \Xsp$. We therefore fix an arbitrary
		$v\in \Xsp$ and compute $(f,v)_{L^2}$ explicitly.
		On each open element $(t_{p},t_{p+1})$, the function $v$ is a polynomial of degree $\le k$, so
		$v^{(k+1)}=0$ there; and across an interior knot its derivatives of order
		$\le\mu$ are continuous, this being the $C^{\mu}$-continuity built into
		$\Xsp$, while its derivatives of order $\mu+1,\dots,k$ may jump. We write
		\begin{equation*}
			[v^{(j)}]_{t_p}:=v^{(j)}(t_p^+)-v^{(j)}(t_p^-), \quad j = \mu+1,\dots,k.
		\end{equation*}
		Let $F:=\J^{k+1}f$ be the antiderivative of $f$, so that $f=F^{(k+1)}$ and $F^{(l)}(0)=0\ \ (l\le k)$. Repeated integration by parts gives
		\[
		(f,v)_{L^2(t_p,t_{p+1})}=\int_{t_{p}}^{t_{p+1}}F^{(k+1)}v
		=\sum_{j=0}^{k}(-1)^{j}\bigl[F^{(k-j)}v^{(j)}\bigr]_{t_{p}}^{t_{p+1}} .
		\]
		Now sum over the $N$ elements. Each interior knot is the right endpoint of
		one element and the left endpoint of the next, so it receives the two
		contributions
		\[
		+\sum_{j=0}^{k}(-1)^{j}F^{(k-j)}(t_{p}^{-})v^{(j)}(t_{p}^{-}),
		\qquad
		-\sum_{j=0}^{k}(-1)^{j}F^{(k-j)}(t_{p}^{+})v^{(j)}(t_{p}^{+}),
		\]
		and since $F^{(k-j)}$ is continuous at $t_{p}$ for $0\le j\le k$, their sum
		is $-\sum_{j=0}^{k}(-1)^{j}F^{(k-j)}(t_{p})[v^{(j)}]_{t_{p}}$. The terms
		at $x=0$ vanish, because $F^{(l)}(0)=0$ for $l\le k$, so
		\begin{equation}\label{eq:IBPstruct}
			(f,v)_{L^{2}}
			=\sum_{j=0}^{k}(-1)^{j}F^{(k-j)}(1)\,v^{(j)}(1)
			-\sum_{p=1}^{N-1}\sum_{j=\mu+1}^{k}(-1)^{j}
			\bigl[v^{(j)}\bigr]_{t_{p}}F^{(k-j)}(t_{p}),
		\end{equation}
		where in the second sum only $j>\mu$ contributes, since
		$[v^{(j)}]_{t_{p}}=0$ for $j\le\mu$. Thus $(f,v)_{L^{2}}$, as a function
		of $v$, is a linear combination of the $k+1+(N-1)\nu=\dim \Xsp$ numbers
		\begin{equation}\label{eq:coordlist}
			v^{(j)}(1)\ (0\le j\le k),
			\qquad
			[v^{(j)}]_{t_{p}}\ (\mu+1\le j\le k,\ 1\le p\le N-1),
		\end{equation}
		with the coefficients $(-1)^{j}F^{(k-j)}(1)$ and
		$-(-1)^{j}F^{(k-j)}(t_{p})$, respectively.

		The map
		\[
		\ell(v):=\Bigl(\bigl(v^{(j)}(1)\bigr)_{j=0}^{k},
		\bigl([v^{(j)}]_{t_p}\bigr)_{\mu+1\le j\le k,\ 1\le p\le N-1}\Bigr)
		\]
		is linear, and its domain and target have the same dimension, so it is an
		isomorphism as soon as it is injective. Injectivity is the implication $\ell(v)=\mathbf 0\Rightarrow
		v=0$, and $\ell(v)=\mathbf 0$ means that all the numbers
		\eqref{eq:coordlist} vanish simultaneously, i.e.
		\begin{equation}\label{eq:lzero}
			v^{(j)}(1)=0\quad(0\le j\le k),\qquad
			[v^{(j)}]_{t_p}=0\quad(\mu+1\le j\le k,\ 1\le p\le N-1).
		\end{equation}
		We deduce $v=0$ from \eqref{eq:lzero} by descending induction on the order of
		differentiation, i.e.\ we show $v^{(j)}=0$ for $j=k,k-1,\dots,0$. Consider $j=k$ first. $v^{(k)}$ is piecewise constant, because $v$ is a
		polynomial of degree $\le k$ on every element. Its jumps $[v^{(k)}]_{t_p}=0$ are listed in \eqref{eq:lzero}, thus vanish, and $v^{(k)}(1)=0$ is also listed there. Hence
		$v^{(k)}=0$. Let $0\le j\le k-1$ and assume $v^{(j+1)}=0$. Then
		$v^{(j)}$ is piecewise constant, and $v^{(j)}(1)=0$ by \eqref{eq:lzero}. Its
		jumps vanish as well: for $j\ge\mu+1$ they are listed in \eqref{eq:lzero},
		and for $j\le\mu$ the function $v^{(j)}$ is continuous by the definition of
		$\Xsp$, so all its jumps are $0$. A piecewise constant
		function with no jumps is a single constant, namely its value at $x=1$,
		which is $0$; hence $v^{(j)}=0$.
		This closes the induction, so $\ell$ is injective, hence an
		isomorphism; in particular the $k+1+(N-1)\nu$ components of $\ell$ form a
		basis of the dual space of $\Xsp$.

		A linear functional on $\Xsp$ vanishes identically if and
		only if all its components relative to the dual basis induced by $\ell$ vanish. Applied
		to \eqref{eq:IBPstruct}, this gives both implications:
		\begin{itemize}[leftmargin=1.6em]
			\item If $F\in\W^{k,\mu}(N)$, then $F^{(k-j)}(1)=0$ for $0\le j\le k$
			and $F^{(k-j)}(t_p)=0$ for $\mu+1\le j\le k$. All coefficients of \eqref{eq:IBPstruct} vanish,
			so $(f,v)_{L^2}=0$ for every $v\in\Xsp$, which means
			$f\in(\Xsp)^{\perp}$.
			\item Conversely, if $(f,v)_{L^2}=0$ for all $v\in\Xsp$, then all
			coefficients of \eqref{eq:IBPstruct} vanish, i.e.\ $F^{(l)}(1)=0$ for $l\le k$ and
			$F^{(l)}(t_p)=0$ for $l\le\nu-1$. Since $F^{(k+1)}=f\in L^2$, we have
			$F\in H^{k+1}$. Thus $F\in\W^{k,\mu}(N)$.
		\end{itemize}
		This completes the proof.
	\end{proof}

	\section{The optimality of the maximally smooth space}\label{sec:floor}

	In this section, we bound $\lamone^{1/2}$ \eqref{eq:Etilde} from above. The
	bound for the maximally smooth space, i.e.
	$E(\Bball,X_{n}^{k,k-1};L^{2})$, asymptotically matches the lower bound
	\eqref{eq:crude} for the Kolmogorov $n$-width. This proves
	Conjecture~\ref{conj:main}, i.e.\ Conjecture 5.1 of \cite{EBH09}, unconditionally.
	We first describe the treatment of the upper bound on $\lamone^{1/2}$ and
	then carry out the detailed analysis.

	By Lemma~\ref{lem:collapse} and Theorem~\ref{thm:structure}, we have
	\begin{equation}\label{eq:lamone}
		\lamone^{1/2}=\sup_{0\ne F\in\Wlat}\frac{\norm{F^{(k+1)}}_{\Hm}}{\norm{F^{(k+1)}}_{L^{2}}},
	\end{equation}
	where the derivative of $F$ is understood in the weak sense. An upper bound for the
	negative norm therefore yields an upper bound for $\lamone^{1/2}$.
	\begin{lemma}\label{lem:hnorm}
		Let $F\in\Wlat$ and $0\le k+1-s\le j\le k+1$. Then
		$\norm{F^{(k+1)}}_{\Hm}\le\norm{F^{(j)}}_{L^{2}}$.
	\end{lemma}

	\begin{proof}
		For $v\in H^{k+1}(0,1)$, a dense subspace of $\Hsp$ because $s\le k+1$, put
		$Q_{i}:=(F^{(i)},v^{(k+1-i)})$; then $(F^{(k+1)},v)_{L^{2}}=Q_{k+1}$, and one
		integration by parts gives $Q_{i}=-Q_{i+1}$ for $i\le k$, the boundary term
		$F^{(i)}v^{(k-i)}\big|_{0}^{1}$ vanishing by the endpoint conditions of $F$. Hence
		$(F^{(k+1)},v)_{L^{2}}=(-1)^{k+1-j}\bigl(F^{(j)},v^{(k+1-j)}\bigr)$, and since
		$k+1-j\le s$,
		\[
		\abs{(F^{(k+1)},v)_{L^{2}}}\le\norm{F^{(j)}}_{L^{2}}\norm{v}_{\Hsp},
		\]
		so $\norm{F^{(k+1)}}_{\Hm}\le\norm{F^{(j)}}_{L^{2}}$ by \eqref{eq:dualmnorm}.
	\end{proof}

	At the level $j=k+1-s$, Lemma~\ref{lem:hnorm} gives
	\begin{equation}\label{eq:reduce}
		\lamone^{1/2}\ \le\
		\sup_{0\ne F\in\Wlat}
		\frac{\norm{F^{(k+1-s)}}_{L^{2}}}{\norm{F^{(k+1)}}_{L^{2}}}.
	\end{equation}
	The error is controlled by the ratio of two derivatives of $F$, and the
	question reduces to the cost of one differentiation on $\Wlat$.
	For $\mu = -1$, this cost is easy to estimate. In this case the conditions
	of Definition~\ref{def:W} make $F^{(l)}(t_p)$ vanish for $l=0,\dots,k$.
	Applying Poincar\'e's inequality on each element $(t_p,t_{p+1})$, one has
	\begin{equation}\label{eq:elementwise-poincare}
		\int_{t_{p}}^{t_{p+1}}\abs{F^{(l)}}^{2}
		\le(N\pi)^{-2}
		\int_{t_{p}}^{t_{p+1}}\abs{F^{(l+1)}}^{2},
		\qquad l=0,\dots,k ,
	\end{equation}
	and summation over the $N$ elements gives
	\[
	\norm{F^{(l)}}_{L^{2}}\le(N\pi)^{-1}\norm{F^{(l+1)}}_{L^{2}},
	\qquad l=0,\dots,k.
	\]
	Hence,
	\begin{equation}\label{eq:elementwise-poincare-2}
		\norm{F^{(k+1-s)}}_{L^{2}}\le(N\pi)^{-s}\norm{F^{(k+1)}}_{L^{2}}.
	\end{equation}
	For a general $\mu>-1$, however, the same computation breaks down. The
	conditions of Definition~\ref{def:W} make $F^{(l)}(t_p)$ vanish for
	$l=0,\dots,\nu-1$ with $\nu=k-\mu$, so the inequality
	\eqref{eq:elementwise-poincare} holds only for $l\le\nu-1$. For the
	maximally smooth space $\nu=1$, this reaches only $l=0$, whereas
	\eqref{eq:elementwise-poincare-2} requires all levels
	$l=k+1-s,\dots,k$. Such a result is nevertheless expected. In what follows
	we prove it from a global rather than an elementwise viewpoint. The
	analysis begins with the expansion below.

	\begin{lemma}\label{lem:sine}
		Let $F\in H^{k+1}(0,1)$ satisfy $F^{(l)}(0)=F^{(l)}(1)=0$ for
		$0\le l\le k$, and put
		\begin{equation}\label{eq:em}
			\psi_m(x):=\sqrt2\sin(m\pi x),\qquad
			e_m:=\int_0^1F(x)\psi_m(x)\,dx\quad(m\ge1).
		\end{equation}
		Then, the sine series
		\begin{equation}\label{eq:sine}
			F=\sum_{m\ge1}e_m\psi_m(x),
		\end{equation}
		converges absolutely and uniformly on $[0,1]$,
		and for $0\le a\le k+1$,
		\begin{equation}\label{eq:norms}
			\|F^{(a)}\|^2_{L^2}=\sum_{m\ge1}(m\pi)^{2a}e_m^2.
		\end{equation}
	\end{lemma}

	\begin{proof}
		Integrating by parts $k+1$ times in the definition \eqref{eq:em} of $e_{m}$, all boundary terms vanishing because $F^{(l)}(0)=F^{(l)}(1)=0$ for $l\le k$, gives
		\[
		e_{m}=(m\pi)^{-(k+1)}c_{m},
		\qquad
		c_{m}:=\int_{0}^{1}F^{(k+1)}(x)\Psi_{m}(x)\,dx,
		\]
		where $\Psi_{m}$ is $\pm\sqrt2\sin(m\pi x)$ or $\pm\sqrt2\cos(m\pi x)$, so
		that $\norm{\Psi_{m}}_{L^{2}}=1$ and the system $\{\Psi_m\}_{m\ge1}$ is
		orthonormal up to signs. Bessel's inequality therefore gives
		$\sum_{m\ge1}c_{m}^{2}\le\norm{F^{(k+1)}}^{2}_{L^{2}}$. Hence, by
		Cauchy--Schwarz,
		\[
		\sum_{m\ge1}\abs{e_{m}}
		\le\Bigl(\sum_{m\ge1}c_{m}^{2}\Bigr)^{1/2}
		\Bigl(\sum_{m\ge1}(m\pi)^{-2(k+1)}\Bigr)^{1/2}
		<\infty,
		\]
		because $2(k+1)\ge2>1$. As $\abs{\psi_{m}}\le\sqrt2$,
		the series $\sum_{m\ge1}e_{m}\psi_{m}$ converges absolutely and uniformly
		on $[0,1]$, and its sum is continuous. That sum is $F$: the series
		converges to $F$ in $L^{2}(0,1)$ because $\{\psi_{m}\}_{m\ge1}$ is a
		complete orthonormal system in $L^{2}(0,1)$, and a continuous function that
		agrees with $F$ in $L^{2}(0,1)$ agrees with it everywhere. This proves
		\eqref{eq:sine}.

		It remains to prove \eqref{eq:norms}. Since the sine series of $F$
		converges uniformly, it may be differentiated term by term in the sense of
		distributions, and for $0\le a\le k+1$ the coefficients of $F^{(a)}$ in
		the orthonormal system
		$\bigl\{\sqrt2\sin(m\pi x+\frac{a\pi}{2})\bigr\}_{m\ge1}$ are
		$(m\pi)^{a}e_{m}$. Up to signs, that system is the sine system when $a$ is even, which is complete in $L^{2}(0,1)$, and the cosine system when $a$ is odd, which together with the constant function is complete; the constant function
		contributes nothing because $F^{(a)}$ has zero mean. Parseval therefore
		gives \eqref{eq:norms}.
	\end{proof}
	By Definition~\ref{def:W}, every $F\in\Wlat$ satisfies the hypotheses of Lemma~\ref{lem:sine}. Evaluating the expansion \eqref{eq:sine} at the interior knots, where
	$F(t_{p})=0$, leads to the following lemma.
	\begin{lemma}\label{lem:alias}
		Let $F\in\Wlat$. Then for every $r=1,\dots,N-1$,
		\begin{equation}\label{eq:alias}
			\sum_{q\ge0}e_{r+2Nq}=\sum_{q\ge0}e_{2N-r+2Nq},
		\end{equation}
		where the coefficients $\{e_m\}_{m\geq1}$ are defined by \eqref{eq:em}.
	\end{lemma}

	\begin{proof}
		The starting point is the vanishing at the knots $F(t_p)=0$ for $p=1,\dots,N-1$, which is the condition of Definition~\ref{def:W} for any $\nu\geq1$, hence for every
		$\mu\le k-1$. By \eqref{eq:sine}, we have
		\begin{equation}\label{eq:knot-values}
			0=F(t_{p})=\sum_{m\ge1}e_m\sqrt2\sin(\frac{m\pi p}{N} ),\qquad p=1,\dots,N-1.
		\end{equation}
		Recall the two elementary identities, for $q\geq 0$ and $r= 1,\dots,N-1 $,
		\begin{gather*}
			\sin((r+2Nq)\pi p/N)=\sin(r\pi p/N),\qquad
			\sin((2N-r+2Nq)\pi p/N)=-\sin(r\pi p/N).
		\end{gather*}
		Regrouping the series \eqref{eq:knot-values} (the terms with $m$ divisible
		by $N$ vanish) gives
		\[
		\sum_{r=1}^{N-1}\sin\Bigl(\frac{r\pi p}{N}\Bigr)
		\Bigl(\sum_{q\ge0}e_{\,r+2Nq}-\sum_{q\ge0}e_{\,2N-r+2Nq}\Bigr)=0,
		\qquad p=1,\dots,N-1 .
		\]
		Clearly, the matrix $\bigl(\sin(r\pi p/N)\bigr)_{r,p=1}^{N-1}$ is invertible, so each bracket vanishes, which is
		\eqref{eq:alias}.
	\end{proof}
	Before presenting the desired result \eqref{eq:elementwise-poincare-2}, we need the following preliminary lemma.
	\begin{lemma}\label{lem:star}
		For $t\in(0,1)$ and every integer $j\ge0$,
		\begin{equation}\label{eq:star}
			\frac{t^{-2j}}{1-t^{2}}\;\ge\;
			\sum_{q\ge1}\Bigl[\frac{(2q-t)^{-2j}}{(2q-t)^{2}-1}
			+\frac{(2q+t)^{-2j}}{(2q+t)^{2}-1}\Bigr],
		\end{equation}
		with equality if and only if $j=0$.
	\end{lemma}

	\begin{proof}
		For $j=0$, partial fractions give
		$\frac{1}{(2q-t)^2-1}=\frac12\bigl(\frac{1}{2q-1-t}-\frac{1}{2q+1-t}\bigr)$,
		and summing over $q\ge1$ telescopes,
		\[
		\sum_{q\ge1}\frac{1}{(2q-t)^2-1}=\frac{1}{2(1-t)}.
		\]
		The same computation gives $\sum_{q\ge1}\frac{1}{(2q+t)^2-1}
		=\frac{1}{2(1+t)}$. Adding, the right-hand side of \eqref{eq:star} equals
		$\frac12\bigl(\frac{1}{1-t}+\frac{1}{1+t}\bigr)=\frac{1}{1-t^2}$, which is
		its left-hand side.

		For $j\ge1$, put $A_{q}:=2q-t$, $w_{q}:=A_{q}^{-2j}$ and
		$c_{q}:=(A_{q}-1)^{-1}$. Since $A_{q+1}-1=A_{q}+1$, the first term of the right-hand side of \eqref{eq:star} is
		\[
		\frac{(2q-t)^{-2j}}{(2q-t)^{2}-1}=
		\frac{A_{q}^{-2j}}{A_{q}^{2}-1}
		=\frac{w_{q}}{(A_{q}-1)(A_{q+1}-1)}
		=\frac12\,w_{q}(c_{q}-c_{q+1}).
		\]
		Abel summation
		\[
		\sum_{q=1}^{Q}w_{q}(c_{q}-c_{q+1})
		=w_{1}c_{1}+\sum_{q=2}^{Q}(w_{q}-w_{q-1})c_{q}-w_{Q}c_{Q+1},
		\]
		with $w_{Q}c_{Q+1}=O(Q^{-2j-1})\to0$ as $Q\to\infty$, gives
		\[
		\sum_{q\ge1}\frac{(2q-t)^{-2j}}{(2q-t)^{2}-1}
		=\frac{w_{1}c_{1}}{2}
		-\frac12\sum_{q\geq2}(w_{q-1}-w_{q})c_{q}
		\ <\frac{w_{1}c_{1}}{2}=\ \frac{(2-t)^{-2j}}{2(1-t)},
		\]
		where every term $(w_{q-1}-w_{q})c_{q}$ is positive. The
		same computation gives
		\begin{equation*}
			\sum_{q\ge1}\frac{(2q+t)^{-2j}}{(2q+t)^{2}-1}
			<\ \frac{(2+t)^{-2j}}{2(1+t)}.
		\end{equation*}
		Adding and writing
		$(2\mp t)^{-2j}=t^{-2j}\bigl(\frac{t}{2\mp t}\bigr)^{2j}$,
		\[
		\sum_{q\ge1}\Bigl[\frac{(2q-t)^{-2j}}{(2q-t)^{2}-1}
		+\frac{(2q+t)^{-2j}}{(2q+t)^{2}-1}\Bigr]
		\ <\ \frac{t^{-2j}}{2}\Bigl[
		\Bigl(\frac{t}{2-t}\Bigr)^{2j}\frac{1}{1-t}
		+\Bigl(\frac{t}{2+t}\Bigr)^{2j}\frac{1}{1+t}\Bigr]
		\ <\ \frac{t^{-2j}}{1-t^{2}},
		\]
		where we use the fact $0<\bigl(\frac{t}{2\mp t}\bigr)^{2j}<1$ when $0<t<1$ and $j\ge1$.
	\end{proof}
	The following lemma will assist us in determining an upper bound on $\lamone^{1/2}$.
	\begin{lemma}\label{lem:class}
		Let $F\in\W^{k,\mu}(N)$ and $0\le j\le k$. Put $h_m:=(m\pi)^je_m$ with $e_m$ defined by \eqref{eq:em}. Then,
		\begin{equation}\label{eq:class}
			\sum_{m\in\mathcal C_i}h_m^2\;\le\;(N\pi)^{-2}
			\sum_{m\in\mathcal C_i}(m\pi)^{2}h_m^2 ,
		\end{equation}
		holds for each of the classes
		\[
		\begin{aligned}
			\mathcal C_r&=\{r+2Nq\}_{q\ge0}\cup\{2N-r+2Nq\}_{q\ge0},\qquad
			1\le r\le N-1,\\
			\mathcal C_0&=\{2N,4N,\dots\},\qquad \mathcal C_N=\{N,3N,5N,\dots\}.
		\end{aligned}
		\]
		Consequently,
		\begin{equation*}
			\|F^{(j)}\|_{L^{2}}\;\le\;(N\pi)^{-1}\|F^{(j+1)}\|_{L^{2}},
		\end{equation*}
		and, for $1\leq s\leq k+1$,
		\begin{equation*}
			\|F^{(k+1-s)}\|_{L^{2}}\le(N\pi)^{-s}\|F^{(k+1)}\|_{L^{2}}.
		\end{equation*}
	\end{lemma}

	\begin{proof}
		For the classes $\mathcal C_0$ and $\mathcal C_N$, $m\ge N$ for
		every member, so $(m\pi)^{2}\ge(N\pi)^{2}$ termwise and \eqref{eq:class} is
		immediate. For a class $\mathcal C_r$, $1\le r\le N-1$, write
		$h:=\{h_m\}_{m\in\mathcal C_r}$ for the vector of coefficients.
		By \eqref{eq:alias}, and since $h_m=(m\pi)^je_m$, the coefficients satisfy
		the single linear relation
		\begin{equation}\label{eq:hyper}
			\sum_{m\in\mathcal C_r}\lambda_mh_m=0,\qquad
			\lambda_m=\epsilon_m(m\pi)^{-j},\qquad
			\epsilon_m=\begin{cases}+1,&m\in\{r+2Nq\}_{q\ge0},\\ -1,&m\in \{2N-r+2Nq\}_{q\ge0}.\end{cases}
		\end{equation}
		Denote $W_m:=(m\pi)^{2}$. Consider the maximum
		\[
		M:=\max\Bigl\{\sum_{m\in\mathcal C_r}h_m^2:\;
		\sum_{m\in\mathcal C_r}W_mh_m^2=1,\;
		\sum_{m\in\mathcal C_r}\lambda_mh_m=0\Bigr\},
		\]
		which is attained. Indeed, since $W_m=(m\pi)^2\to\infty$, the closed
		unit ball of the weighted space
		$\{h:\sum_mW_mh_m^2\le1\}$ is compact in $\ell^2$: for every $m^*$,
		\[
		\sum_{m>m^*}h_m^2\le W_{m^*+1}^{-1}\sum_{m>m^*}W_mh_m^2\le W_{m^*+1}^{-1},
		\]
		so its tail is uniformly small. Replacing the sphere by this ball
		therefore gives a compact admissible set, and $\sum_mh_m^2$ is
		continuous on it, so the maximum is attained. A maximizer lies on the
		sphere, since otherwise rescaling by $t>1$ would increase it. The
		gradients of the two constraints are the vectors
		$(2W_mh_m)_{m\in\mathcal C_r}$ and $(\lambda_m)_{m\in\mathcal C_r}$, and
		they are linearly independent, that is, there is no single constant $c$
		with $2W_mh_m=c\lambda_m$ for every $m$. Indeed, such a $c$ would give,
		on multiplying by $\lambda_m/W_m$, summing, and recalling the constraint
		$\sum_{m\in\mathcal C_r}\lambda_mh_m=0$,
		\[
			0=2\sum_{m\in\mathcal C_r}\lambda_mh_m
			=c\sum_{m\in\mathcal C_r}\frac{\lambda_m^2}{W_m},
		\]
		hence $c=0$. Every term $\lambda_m^2/W_m$ is strictly positive, since
		$\lambda_m=\pm(m\pi)^{-j}\ne0$ and $0<W_m=(m\pi)^2<\infty$ for each $m$,
		so the sum cannot vanish. But then $h=0$, which contradicts the
		normalization $\sum_mW_mh_m^2=1$. The multiplier rule therefore applies.
		The maximizer satisfies, with multipliers $\alpha,\beta$,
		\[
		2h_m-2\alpha W_mh_m-\beta\lambda_m=0,
		\]
		hence
		\begin{equation}\label{eq:lagrange}
			h_m=\frac{\beta\lambda_m}{2(1-\alpha W_m)} .
		\end{equation}
		If $\beta=0$, then $h_m(1-\alpha W_m)=0$ for all $m$. Since the values $W_m$ are distinct, $h$ is supported on at most one index; the constraint $\sum\lambda_mh_m=0$ then forces $h=0$,
		which contradicts the normalization $\sum_mW_mh_m^2=1$. Hence $\beta\ne0$, and
		inserting \eqref{eq:lagrange} into \eqref{eq:hyper} and dividing by
		$\beta/2$,
		\begin{equation}\label{eq:Phi}
			\sum_{m\in\mathcal C_r}\frac{\lambda_m^2}{1-\alpha W_m}=0,
			\qquad\text{i.e.}\qquad
			\Phi(\alpha):=\sum_{m\in\mathcal C_r}\frac{(m\pi)^{-2j}}{1-\alpha(m\pi)^{2}}=0.
		\end{equation}
		For such a critical point,
		$\sum_m\bigl(1-\alpha W_m\bigr)h_m^2=\frac{\beta}{2}\sum_m\lambda_mh_m=0$ by
		\eqref{eq:lagrange} and \eqref{eq:hyper}, so
		\[
		\sum_{m\in\mathcal C_r}h_m^2=\alpha\sum_{m\in\mathcal C_r}W_mh_m^2=\alpha .
		\]
		Consequently, the value of the maximum is a root of \eqref{eq:Phi}. Conversely,
		every root of \eqref{eq:Phi} comes from an admissible point of value $\alpha$:
		let $\alpha$ be a root and fix any $\beta\ne0$. As a root of $\Phi$, $\alpha$ is
		not one of its poles, that is, $1-\alpha W_m\ne0$ for every $m$, so
		\[
		h_m:=\frac{\beta\lambda_m}{2(1-\alpha W_m)}
		\]
		is well defined, and $h\ne0$ because $\lambda_m\ne0$. By construction,
		$(1-\alpha W_m)h_m=\frac{\beta\lambda_m}{2}$, hence
		\[
		\sum_m(1-\alpha W_m)h_m^2=\frac{\beta}{2}\sum_m\lambda_mh_m
		=\frac{\beta^2}{4}\Phi(\alpha)=0 ,
		\]
		the last equality because $\alpha$ is a root of $\Phi$. Therefore
		$\sum_mh_m^2=\alpha\sum_mW_mh_m^2$, and $h$ may be rescaled so that
		$\sum_mW_mh_m^2=1$; both this relation and the constraint \eqref{eq:hyper} are
		homogeneous, so they survive the rescaling, and the value becomes
		$\sum_mh_m^2=\alpha$. Hence every root of \eqref{eq:Phi} is at most $M$, while
		$M$ itself is a root; that is, $M=\alpha_\ast$, the largest root.

		The poles of $\Phi$ are the points $\alpha=(m\pi)^{-2}$ with $m\in\mathcal C_r$; order
		them decreasingly. On each interval between consecutive poles,
		\[
		\Phi'(\alpha)=\sum_{m\in\mathcal C_r}
		\frac{(m\pi)^{2-2j}}{\bigl(1-\alpha(m\pi)^2\bigr)^2}\ >\ 0 ,
		\]
		so $\Phi$ is strictly increasing there, while at the pole $\alpha=(m\pi)^{-2}$,
		the summand with index $m$ tends to $+\infty$ from the left and to $-\infty$ from the
		right, all other summands staying bounded there. Hence on each such interval,
		$\Phi$ runs from $-\infty$ up to $+\infty$ and has exactly one zero in it. The
		two largest poles come from the two smallest elements of $\mathcal C_r$, so the largest zero $\alpha_\ast$ lies in
		\[
		\Bigl(\frac{1}{((2N-r)\pi)^2},\ \frac{1}{(r\pi)^2}\Bigr).
		\] The point $\alpha_0:=(N\pi)^{-2}$ lies in the
		same interval, because $r<N<2N-r$; and since $\Phi$ increases there and
		$\Phi(\alpha_\ast)=0$,
		\begin{equation}\label{eq:Phi-sign}
			\Phi(\alpha_0)\ge0\quad\Longrightarrow\quad\alpha_0\ge\alpha_\ast=M
			\quad\Longrightarrow\quad M\le(N\pi)^{-2}.
		\end{equation}
		Therefore, it remains to show $ \Phi(\alpha_0)\ge0$.
		Multiplying $\Phi(\alpha_0)$ by $(N\pi)^{2j}$ and writing $m=Nu$ turns
		the summand $\frac{(m\pi)^{-2j}}{1-\alpha_0(m\pi)^{2}}$ into
		$\frac{u^{-2j}}{1-u^{2}}$, and the class $\mathcal C_r$ into the numbers
		$u=t+2q$ with $q\ge0$ and $u=2q-t$ with $q\ge1$, where $t=r/N$; hence
		\[
		(N\pi)^{2j}\Phi(\alpha_0)
		=\sum_{q\ge0}\frac{(t+2q)^{-2j}}{1-(t+2q)^{2}}
		+\sum_{q\ge1}\frac{(2q-t)^{-2j}}{1-(2q-t)^{2}} .
		\]
		The term $q=0$ of the first sum is $\frac{t^{-2j}}{1-t^{2}}$, so
		\[
		(N\pi)^{2j}\Phi(\alpha_0)
		=\frac{t^{-2j}}{1-t^{2}}
		-\sum_{q\ge1}\Bigl[\frac{(2q-t)^{-2j}}{(2q-t)^{2}-1}
		+\frac{(2q+t)^{-2j}}{(2q+t)^{2}-1}\Bigr]\;\ge\;0
		\]
		by Lemma~\ref{lem:star}. Thus, $\Phi(\alpha_0)\ge0$, and $M\le(N\pi)^{-2}$ by
		\eqref{eq:Phi-sign}; that is, for the normalized $h$,
		$\sum_{\mathcal C_r}h_m^2\le(N\pi)^{-2}\sum_{\mathcal C_r}W_mh_m^2$.
		For an arbitrary admissible $h\ne0$, put
		$\sigma:=\sum_{\mathcal C_r}W_mh_m^2>0$ and $\tilde h:=h/\sqrt{\sigma}$. Then
		$\tilde h$ is admissible, because the constraint is homogeneous, and
		$\sum_{\mathcal C_r}W_m\tilde h_m^2=1$, so the inequality just proved applies
		to $\tilde h$; after multiplying by $\sigma$ it reads
		$\sum_{\mathcal C_r}h_m^2\le(N\pi)^{-2}\sum_{\mathcal C_r}W_mh_m^2$, which is
		\eqref{eq:class} for the class $\mathcal C_r$. Together with the classes
		$\mathcal C_0,\mathcal C_N$ treated at the beginning, this proves
		\eqref{eq:class} for every class.
	\end{proof}

	Two features of Lemma~\ref{lem:class} together with its proof explain both its strength and its limitation. First, only the vanishing $F(t_{p})=0$ at the knots was used, not the
	vanishing of the higher derivatives in Definition~\ref{def:W}, $F^{(l)}(t_{p})=0$,
	$1\le l\le\nu-1$. Hence the lemma holds for the whole family
	$\Wlat$, uniformly in $\mu$. Second, the discarded conditions can only make the true maximum smaller:
	the vanishing of the higher derivatives may add constraints to the maximum above, and adding constraints can only lower a maximum. Therefore, this upper bound of the maximum may be much looser when $\mu<k-1$.

	At the end of this section, we show the optimality of the maximally smooth space for the Kolmogorov $n$-width.
	\begin{theorem}\label{thm:main}
			Let $k\ge s-1$ and $1\leq \nu\leq k+1$ with $\nu=k-\mu$; let $\Xsp$ be
		the spline space \eqref{eq:X}, of dimension $n=N\nu+\mu+1$. Then the ratio $ \Lam$ given by \eqref{eq:introLam} has the following upper bound
		\begin{equation}\label{eq:main}
			\Lam\ \le\ \nu^{s}\Bigl(1+\frac{s(\mu+1)}{n}+O(n^{-2})\Bigr).
		\end{equation}
		In particular, for the maximally smooth space ($\nu=1$),
		\begin{equation}\label{eq:mainsmooth}
			\Lams\ \le\ 1+\frac{sk}{n}+O(n^{-2}),
		\end{equation}
		and Conjecture~\ref{conj:main} holds: $\Lams\to1$ as $n \to\infty$.
	\end{theorem}

	\begin{proof}
		By Lemma~\ref{lem:class} and \eqref{eq:reduce},
		$\lamone^{1/2}\le(N\pi)^{-s}$; by \eqref{eq:crude},
		$\dn\ge(n\pi)^{-s}(1-O(n^{-2}))$. Dividing the two and using
		$n=N\nu+\mu+1$ gives \eqref{eq:main}. For $\nu=1$, this reads
		\eqref{eq:mainsmooth}, and $\Lam\ge1$ yields $\Lams\to1$.
	\end{proof}

	Although the upper bound \eqref{eq:main} is not sharp for all spline spaces $\Xsp$, it already shows that the maximally smooth space is asymptotically optimal: its error asymptotically matches the Kolmogorov $n$-width as $n \to\infty$.

	\section{Lower bounds for the ratio \texorpdfstring{$\Lam$}{Lambda}}\label{sec:nu}
	In this section, we derive lower bounds for the error $\lamone^{1/2}$ of
	the spline space, and hence for the ratio $\Lam$. By the identity
	\eqref{eq:lamone}, a single function suffices: if $0\ne F_{0}\in\Wlat$, then
	\begin{equation}\label{eq:onelb}
		\lamone^{1/2}\ \ge\ \frac{\norm{F_{0}^{(k+1)}}_{\Hm}}
		{\norm{F_{0}^{(k+1)}}_{L^{2}}} .
	\end{equation}
	Definition~\ref{def:W} specifies the conditions that $F_{0}$ must satisfy. Its
	derivatives up to order $\nu-1$ have to vanish at the interior knots
	$t_{p}=p/N$, $1\le p\le N-1$, and its derivatives up to order $k$ have to
	vanish at the endpoints $0$ and $1$. This suggests a natural choice. The
	factor $\sin^{\nu}(MN\pi x)$, with $M\ge1$ an integer, makes the
	derivatives vanish at the interior knots, and it also gives
	$F_{0}^{(l)}(x)|_{x\in\{0,1\}}=0$ for $0\le l\le\nu-1$, so the endpoint
	conditions of order $\nu,\dots,k$ still require a second factor. We take
	the second factor as $b(x)$, with $b\in C^{k+1}([0,1])$. By Leibniz, it
	suffices that it satisfy
	\begin{equation}\label{eq:def-b}
		b^{(l)}(0)=b^{(l)}(1)=0,\quad 0\le l\le k-\nu=\mu.
	\end{equation}
	In what follows, we therefore consider
	\begin{equation}\label{eq:envclass}
		F_{0}(x)=\sin^{\nu}(\omega x)\,b(x),
		\qquad \omega:=MN\pi,
	\end{equation}
	with $M\ge1$ an integer and $b\in C^{k+1}([0,1])$ satisfying the stated
	conditions \eqref{eq:def-b}. Choosing two test functions for the $\Hm$-norm,
	the single function $\delta=F_{0}^{(k+1)}$ and the cosine sequence, we
	analyze the corresponding lower bounds in the two subsections below. The
	numerical experiments at the end of the paper show that the bound obtained
	from this choice of $F_{0}$ is sharp whenever the two lower bounds coincide.

	\subsection{A crude lower bound for the ratio}\label{sec:sharp}
	We begin with the expansion of the oscillating factor $\sin^{\nu}\theta$,
	where $\theta:=\omega x$. Since $\sin\theta=(e^{i\theta}-e^{-i\theta})/(2i)$,
	we have
	\begin{equation}\label{eq:pulse}
		\sin^{\nu}\theta
		=\Bigl(\frac{e^{i\theta}-e^{-i\theta}}{2i}\Bigr)^{\nu}
		=\sum_{m=0}^{\nu}\binom{\nu}{m}\frac{(-1)^{m}}{(2i)^{\nu}}
		e^{i(\nu-2m)\theta}
		:= \sum_{j\in J_{\nu}}c_{\nu,j}e^{ij\theta},
	\end{equation}
	with the corresponding coefficients $c_{\nu,j}$ and $J_{\nu}=\{\nu,\nu-2,\nu-4,\cdots,-\nu \}$.
	Define
	\begin{equation}\label{eq:Abump}
		A_{\nu,r}:=\sum_{j\in J_{\nu}}|c_{\nu,j}|^2j^{2r}
		=4^{-\nu}\sum_{m=0}^{\nu}\binom{\nu}{m}^2|\nu-2m|^{2r},\qquad r\geq 0.
	\end{equation}
	Here and below we use the convention $0^0=1$; this affects only the case
	$r=0$.
	\begin{lemma}\label{lem:moments}
		Let $1\le\nu\le k+1$. Take $F_{0}$ as in \eqref{eq:envclass}, where $b(x)$ is assumed further to be sufficiently smooth.
		Then,
		\begin{equation}\label{eq:momconcl}
			\|F_{0}^{(r)}\|^2_{L^2}
			=\omega^{2r}\|b\|^2 _{L^{2}}A_{\nu,r}\bigl(1+O(\omega^{-2})\bigr),
		\end{equation}
		where $A_{\nu,r}$ is defined by \eqref{eq:Abump}.
	\end{lemma}

	\begin{proof}
		Write $\theta=\omega x$ and expand by Leibniz,
		\begin{equation}\label{eq:leib}
			F_{0}^{(r)}=\sum_{l=0}^{r}\binom rl\omega^l\,b^{(r-l)}(x)\,
			\Bigl(\frac{d^l}{d\theta^l}\sin^{\nu}\Bigr)(\omega x),
		\end{equation}
		where the factor $\bigl(\frac{d^l}{d\theta^l}\sin^{\nu}\bigr)(\omega x)$ is
		the $l$-th derivative of $\sin^{\nu}(\theta)$ with respect to $\theta$,
		evaluated at $\theta=\omega x$. Thus
		$\|F_{0}^{(r)}\|^2_{L^2}$ is the sum over $l,l'$ of
		\begin{equation}\label{eq:leibterm}
			\binom rl\binom r{l'}\omega^{l+l'}\int_0^1b^{(r-l)}(x)\,b^{(r-l')}(x)\,
			\Bigl(\frac{d^l}{d\theta^l}\sin^{\nu}\,\cdot\,
			\frac{d^{l'}}{d\theta^{l'}}\sin^{\nu}\Bigr)(\omega x)\,dx .
		\end{equation}
		Here $g:=b^{(r-l)}b^{(r-l')}$, and
		$(\frac{d^l}{d\theta^l}\sin^{\nu}\,\cdot\,
		\frac{d^{l'}}{d\theta^{l'}}\sin^{\nu})(\omega x)$ is a linear combination
		of sines and cosines. Let $g\in C^{2}([0,1])$. For $\ell\ne0$, we have
		\[
		\int_0^1g(x)\cos(2\ell\omega x)\,dx
		=\frac{g'(1)-g'(0)}{(2\ell\omega)^2}+O(\omega^{-2})
		=O(\omega^{-2}),
		\]
		since $\sin(2\ell\omega)=0$ removes the boundary term of the first
		integration by parts, and
		\[
		\int_0^1g(x)\sin(2\ell\omega x)\,dx
		=-\frac{g(1)\cos(2\ell\omega)-g(0)}{2\ell\omega}+O(\omega^{-1})
		=O(\omega^{-1}).
		\]

		We consider the parity of
		$\bigl(\frac{d^{l}}{d\theta^{l}}\sin^{\nu}\cdot
		\frac{d^{l'}}{d\theta^{l'}}\sin^{\nu}\bigr)(\omega x)$. Since
		$\sin^{\nu}(-\theta)=(-1)^{\nu}\sin^{\nu}(\theta)$, differentiating this
		identity $l$ times, the chain rule contributing $(-1)^{l}$, gives
		\[
		\Bigl(\frac{d^{l}}{d\theta^{l}}\sin^{\nu}\Bigr)(-\theta)
		=(-1)^{\nu+l}\Bigl(\frac{d^{l}}{d\theta^{l}}\sin^{\nu}\Bigr)(\theta),
		\]
		so the two factors of the product have the parities $\nu+l$ and
		$\nu+l'$: the product is an even function of $\theta$ when $|l-l'|$ is
		even and an odd one when $|l-l'|$ is odd. An even trigonometric polynomial
		contains only the constant and cosines, while an odd one contains only
		sines. Hence
		\[
		\Bigl(\frac{d^{l}}{d\theta^{l}}\sin^{\nu}\,\cdot\,
		\frac{d^{l'}}{d\theta^{l'}}\sin^{\nu}\Bigr)(\omega x)=
		\begin{cases}
			\displaystyle d_0(l,l')+\sum_{1\le\abs\ell\le\nu}
			\alpha_{\ell}\cos(2\ell\omega x), & |l-l'|\ \text{even},\\[1.2mm]
			\displaystyle \sum_{1\le\abs\ell\le\nu}\beta_{\ell}
			\sin(2\ell\omega x), & |l-l'|\ \text{odd},
		\end{cases}
		\]
		where $d_0(l,l')\in\R$ is the constant term. By \eqref{eq:pulse},
		\[
		d_0(l,l'):=\sum_{\substack{j,j'\in J_{\nu}\\ j+j'=0}}
		c_{\nu,j}c_{\nu,j'}(ij)^{l}(ij')^{l'} .
		\]
		With this form, the terms \eqref{eq:leibterm} are estimated by the two
		integrals above. If $|l-l'|$ is odd, every term is a combination of
		$\int_0^1g(x)\sin(2\ell\omega x)\,dx$; if $|l-l'|$ is even, every term
		is a combination of $\int_0^1g(x)\cos(2\ell\omega x)\,dx$ and
		$d_0(l,l')\int_0^1g$.
		Hence, for every $g\in C^{2}([0,1])$,
		\begin{equation}\label{eq:crossterm}
			\int_0^1g(x)\,
			\Bigl(\frac{d^{l}}{d\theta^{l}}\sin^{\nu}\,\cdot\,
			\frac{d^{l'}}{d\theta^{l'}}\sin^{\nu}\Bigr)(\omega x)\,dx
			=\begin{cases}
				\displaystyle d_0(l,l')\int_0^1g+O(\omega^{-2}), & |l-l'|\ \text{even},\\[1.2mm]
				O(\omega^{-1}), & |l-l'|\ \text{odd}.
			\end{cases}
		\end{equation}

		Now we sum \eqref{eq:leibterm} over the pairs $l,l'$. The
		diagonal term $l=l'=r$ of \eqref{eq:leib} contributes
		$\omega^{2r}\norm{b}^{2}_{L^{2}}d_0(r,r)(1+O(\omega^{-2}))$, by \eqref{eq:crossterm}
		with $g=b^{2}$. Every off-diagonal term contributes $O(\omega^{2r-2})$.
		If $l-l'$ is even, the $d_0$-part carries the weight
		$\omega^{l+l'}\le\omega^{2r-2}$. If $l-l'$ is odd, the whole term is
		$O(\omega^{l+l'-1})$ with $l+l'\le2r-1$. Thus,
		\begin{equation*}
			\norm{F_{0}^{(r)}}^{2}_{L^{2}}
			=\omega^{2r}\norm{b}^{2}_{L^{2}}d_0(r,r)\bigl(1+O(\omega^{-2})\bigr),
		\end{equation*}
		and it remains to recognize the factor $d_0(r,r)$. On the diagonal,
		$d_0(l,l)$ pairs $c_{\nu,j}$, $c_{\nu,j'}$ and $(ij)^{l}$, $(ij')^{l}$ with $j'=-j$, where
		$(ij)^{l}(-ij)^{l}=j^{2l}$; since $\sin^{\nu}$ is
		real, its coefficients satisfy $c_{\nu,-j}=\overline{c_{\nu,j}}$, so that
		$c_{\nu,j}c_{\nu,-j}=\abs{c_{\nu,j}}^{2}$. Hence
		\begin{equation*}
			d_0(l,l)=\sum_{j\in J_{\nu}}\abs{c_{\nu,j}}^{2}j^{2l}=A_{\nu,l}
		\end{equation*}
		by the definition \eqref{eq:Abump}, and in particular $d_0(r,r)=A_{\nu,r}$.
		Inserting this above gives \eqref{eq:momconcl}.
	\end{proof}

	This yields a lower bound for $\lamone^{1/2}$.
	\begin{lemma}\label{lem:crude}
		Let $1\le\nu\le k+1$. Then,
		\begin{equation}\label{eq:cs}
			\lamone^{1/2}\ \ge \ (N\pi)^{-s}\Bigl(\frac{\Ab{k+1}}{\Ab{k+1+s}}\Bigr)^{1/2}
			\bigl(1-O(N^{-2})\bigr),
		\end{equation}
		where $A_{\nu,r}$ is defined by \eqref{eq:Abump}.
	\end{lemma}

	\begin{proof}
		By \eqref{eq:onelb},
		\[
		\lamone^{1/2}\ \ge\
		\frac{\norm{F_{0}^{(k+1)}}_{\Hm}}{\norm{F_{0}^{(k+1)}}_{L^{2}}}
		\ \ge\
		\frac{\norm{F_{0}^{(k+1)}}_{L^{2}}}{\norm{F_{0}^{(k+1)}}_{H^{s}}},
		\]
		where the second inequality holds because $\norm{\cdot}_{\Hm}$ defined by \eqref{eq:dualmnorm} is a supremum over
		$0\ne\delta\in\Hsp$, and the choice $\delta=F_{0}^{(k+1)}$ there gives
		\[
		\norm{F_{0}^{(k+1)}}_{\Hm}\ \ge\
		\frac{\norm{F_{0}^{(k+1)}}^{2}_{L^{2}}}{\norm{F_{0}^{(k+1)}}_{H^{s}}} .
		\]
		Take now $F_{0}$ to be \eqref{eq:envclass}, with, for instance,
		$b(x)=x^{\mu+1}(1-x)^{\mu+1}$, which satisfies \eqref{eq:def-b}. By Lemma~\ref{lem:moments},
		\[
		\frac{\norm{F_{0}^{(k+1)}}^{2}_{L^{2}}}
		{\norm{F_{0}^{(k+1)}}^{2}_{H^{s}}}
		=\frac{\Ab{k+1}}{\sum_{a=0}^{s}\omega^{2a}\Ab{k+1+a}}
		\bigl(1+O(\omega^{-2})\bigr),
		\]
		hence, as a lower bound,
		\[
		\lamone^{1/2}\ \ge\
		\Bigl(\frac{\Ab{k+1}}{\sum_{a=0}^{s}\omega^{2a}\Ab{k+1+a}}\Bigr)^{1/2}
		\bigl(1-O(\omega^{-2})\bigr).
		\]
		In the denominator of this ratio, the dominant term is the top level $a=s$:
		the sum equals $\omega^{2s}\Ab{k+1+s}\bigl(1+O(\omega^{-2})\bigr)$, the lower
		levels $a<s$ entering the remainder only. Hence
		\[
		\lamone^{1/2}\ \ge\ \omega^{-s}
		\Bigl(\frac{\Ab{k+1}}{\Ab{k+1+s}}\Bigr)^{1/2}
		\bigl(1-O(\omega^{-2})\bigr),
		\]
		and this bound decreases with $\omega$, so it is largest at the smallest $M$, that is, at $\omega=N\pi$:
		\[
		\lamone^{1/2}\ \ge\ (N\pi)^{-s}
		\Bigl(\frac{\Ab{k+1}}{\Ab{k+1+s}}\Bigr)^{1/2}
		\bigl(1-O(N^{-2})\bigr),
		\]
		which is \eqref{eq:cs}.
	\end{proof}

	For $\nu=1$, the set $J_{\nu}$ in \eqref{eq:Abump} is $\{-1,1\}$.
	For every $r$, $A_{1,r}=4^{-1}\bigl(1^{2r}+(-1)^{2r}\bigr)=1/2$, and hence
	\begin{equation*}
		\frac{A_{1,k+1}}{A_{1,k+1+s}}\ =\ 1.
	\end{equation*}
	So the lower bound \eqref{eq:cs} for the maximally
	smooth space reads
	\begin{equation*}
		\lamone^{1/2}\ge(N\pi)^{-s}\bigl(1-O(N^{-2})\bigr).
	\end{equation*}
	Lemma~\ref{lem:class} together with \eqref{eq:reduce}
	gives the matching upper bound $\lamone^{1/2}\le(N\pi)^{-s}$. The two leading
	terms coincide, so $(N\pi)^{-s}$ is sharp. We record this as a corollary.

	\begin{corollary}\label{thm:nu1}
		Let $k\ge s-1$. For the maximally smooth space $\Kh$,
		\begin{equation*}
			\lamone^{1/2}=E(\Bball,\Kh;L^{2})=(N\pi)^{-s}\bigl(1+O(N^{-2})\bigr).
		\end{equation*}
	\end{corollary}

	The corollary above gives the sharp form of the error $\lamone^{1/2}$;
	combined with the conjectured sharp form of the Kolmogorov $n$-width
	(Conjecture~\ref{conj:shift}), it gives the corresponding sharp conjecture for
	the ratio $\Lams$ itself.

	\begin{conjecture}\label{conj:ratio}
		Let $k\ge s-1$. Then,
		\begin{equation*}
			\Lams=1+\frac{s(2k+1-s)}{2n}+O(n^{-2}).
		\end{equation*}
		Here, the coefficient $s(2k+1-s)/2$ is strictly smaller than the $sk$ of
		\eqref{eq:mainsmooth} as soon as $s>1$, with equality only for $s=1$.
	\end{conjecture}

	By Theorem~\ref{thm:main} and $\Lam\ge1$, the ratio $\Lams$ of the maximally
	smooth space is $1+O(n^{-1})$. The
	Kolmogorov $n$-width reads
	\begin{equation}\label{eq:dnform}
		\dn=(n\pi)^{-s}\bigl(1+O(n^{-1})\bigr).
	\end{equation}
	Dividing the lower bound \eqref{eq:cs} of Lemma~\ref{lem:crude} by
	\eqref{eq:dnform}, and using
	$n=N\nu+\mu+1=\nu N\bigl(1+O(n^{-1})\bigr)$, turns it into a lower bound for the ratio $\Lam$.

	\begin{theorem}\label{thm:nucrude}
		Let $1\le\nu\le k+1$ and $k\ge s-1$. Then,
		\begin{equation}\label{eq:nucrude}
			\Lam\ \ge\ \nu^{s}\Bigl(\frac{\Ab{k+1}}{\Ab{k+1+s}}\Bigr)^{1/2}
			\bigl(1-O(n^{-1})\bigr),
		\end{equation}
		where $A_{\nu,r}$ is defined by \eqref{eq:Abump}.
	\end{theorem}

	The bound \eqref{eq:nucrude} is strictly larger than $1$ only for $\nu>2$. At
	$\nu=2$, the set $J_{\nu}$ in \eqref{eq:Abump} is $\{-2,0,2\}$. For every
	$r\ge1$, $A_{2,r}=4^{-2}\bigl(2^{2r}+(-2)^{2r}\bigr)=2^{\,2r-3}$, and hence
	\begin{equation*}
		\nu^{s}\Bigl(\frac{\Ab{k+1}}{\Ab{k+1+s}}\Bigr)^{1/2}=2^{s}2^{-s}=1.
	\end{equation*}
	For $\nu>2$, the lower bound \eqref{eq:nucrude} is strictly larger than $1$.
	Indeed, by \eqref{eq:Abump},
	\begin{equation*}
		\nu^{2s}\Ab{k+1}-\Ab{k+1+s}
		=4^{-\nu}\sum_{m=0}^{\nu}\binom{\nu}{m}^{2}
		|\nu-2m|^{2(k+1)}\bigl(\nu^{2s}-|\nu-2m|^{2s}\bigr),
	\end{equation*}
	where every term on the right is non-negative. The term
	$m=1$,
	\[
	\nu^{2}(\nu-2)^{2(k+1)}\bigl(\nu^{2s}-(\nu-2)^{2s}\bigr),
	\]
	is positive since $\nu>2$. Hence
	$\nu^{2s}\Ab{k+1}>\Ab{k+1+s}$, which is the strict inequality
	$\nu^{s}\bigl(\Ab{k+1}/\Ab{k+1+s}\bigr)^{1/2}>1$ of the bound. This comparison reveals the superiority of the maximally smooth spline space over the less smooth spline spaces ($\nu>2$) in the sense of $n$-widths. A finer comparison, in which the $\Hm$-norm is tested against more
	functions than the single $\delta=F_{0}^{(k+1)}$ used so far, is the
	subject of the next subsection.

	\subsection{A sharper lower bound for the ratio}\label{sec:nu-sharp}

	As announced at the end of the previous subsection, the norm
	$\norm{F_{0}^{(k+1)}}_{\Hm}$ is now tested against more functions than the
	single $\delta=F_{0}^{(k+1)}$ used there. The supremum \eqref{eq:dualmnorm}
	defining that norm runs over all of $\Hsp$, and the more functions it is
	tested against, the sharper the lower bound it yields; in this subsection we
	test it against the cosines and evaluate the resulting series.

	We work with the cosines, normalized in $L^{2}$,
	\begin{equation*}
		\varphi_{m}:=\begin{cases}1,&m=0,\\ \sqrt2\cos(m\pi x),&m\ge1,\end{cases}
		\qquad
		\Ss(m):=\sum_{a=0}^{s}(m\pi)^{2a},
	\end{equation*}
	a system orthonormal in $L^{2}(0,1)$ and $H^{s}$-orthogonal with
	$\norm{\varphi_{m}}^{2}_{H^{s}}=\Ss(m)$. The derivatives of $\varphi_{m}$
	are sines and cosines of the frequency $m\pi$, which is the computation
	already made in the proof of Lemma~\ref{lem:minmax}. A cosine combination
	$v=\sum_{m\le M}d_{m}\varphi_{m}$ therefore has
	$\norm{v}^{2}_{H^{s}}=\sum_{m\le M}\Ss(m)d_{m}^{2}$ and
	$(g,v)_{L^{2}}=\sum_{m\le M}d_{m}(g,\varphi_{m})$. Cauchy--Schwarz gives
	\[
	\frac{(g,v)_{L^{2}}}{\norm{v}_{H^{s}}}
	=\frac{1}{\norm{v}_{H^{s}}}\sum_{m\le M}\bigl(\sqrt{\Ss(m)}\,d_{m}\bigr)
	\Bigl(\frac{(g,\varphi_{m})}{\sqrt{\Ss(m)}}\Bigr)
	\le \Bigl(\sum_{m\le M}\frac{(g,\varphi_{m})^{2}}{\Ss(m)}\Bigr)^{1/2},
	\]
	with equality for $d_{m}=(g,\varphi_{m})/\Ss(m)$. The value on the right is
	therefore attained by a cosine combination, so the supremum
	\eqref{eq:dualmnorm}, which runs over all of $\Hsp$, is at least that large
	for every $M$. The partial sums increase, and letting $M\to\infty$ gives
	\begin{equation}\label{eq:cosnorm}
		\norm{g}^{2}_{\Hm}\ \ge\ \sum_{m\ge0}\frac{(g,\varphi_{m})^{2}}{\Ss(m)} .
	\end{equation}
	For $s=1$, one integration by parts gives
	$(u,\varphi_{m})_{H^{1}}=\Ss(m)(u,\varphi_{m})_{L^{2}}$ for every
	$u\in H^{1}$, and the cosine span is dense in $H^1$ as well. Hence the
	supremum \eqref{eq:dualmnorm} equals the supremum over cosine combinations,
	which the estimate above evaluates as the series on the right of
	\eqref{eq:cosnorm}. Therefore \eqref{eq:cosnorm} is an identity.
	For $s\ge2$, this identity need not hold. Every cosine combination has vanishing first
	derivatives at the endpoints, and since $\Hsp\subset C^{1}$ for $s\ge2$ the
	values $v'(0)$ and $v'(1)$ are continuous on $\Hsp$. The condition survives
	the limit, and the closure of the cosine span lies in the proper subspace
	$\{v\in\Hsp:\ v'(0)=v'(1)=0\}$. The cosine span is therefore not dense in
	$\Hsp$, and the supremum may be strictly larger. The argument below uses
	\eqref{eq:cosnorm} only as a lower bound for the dual norm, and the series
	on its right is the object of the evaluation.
	Since
	$\varphi_{m}^{(k+1)}=(m\pi)^{k+1}\Psi_{m}$ for $m\ge1$ with
	$\Psi_{m}=\sqrt2\cos\bigl(m\pi x+(k+1)\pi/2\bigr)$, $k+1$ integrations by
	parts give
	\begin{equation*}
		(F_{0}^{(k+1)},\varphi_{m})_{L^{2}}
		=(-1)^{k+1}(m\pi)^{k+1}c_{m},
		\qquad
		c_{m}:=\int_{0}^{1}F_{0}\,\Psi_{m},
	\end{equation*}
	all boundary terms vanishing because $F_0$ satisfies the endpoint conditions, which follow from \eqref{eq:def-b} and the sine factor; and
	$(F_{0}^{(k+1)},\varphi_{0})=\int_{0}^{1}F_{0}^{(k+1)}
	=F_{0}^{(k)}(1)-F_{0}^{(k)}(0)=0$, so the sum over $m$ starts at $m=1$.
	Hence, by
	\eqref{eq:cosnorm},
	\begin{equation}\label{eq:cosnormE}
		\norm{F_{0}^{(k+1)}}^{2}_{\Hm}
		\ \ge\ \sum_{m\ge1}\frac{(m\pi)^{2(k+1)}}{S_s(m)}c_{m}^{2},
	\end{equation}
	and the problem is reduced to the coefficients $c_{m}$ of $F_{0}$.

	To compute them we expand $F_{0}$. From here on, we take the smallest
	value of $M$, that is, $\omega=N\pi$ in \eqref{eq:envclass}: the evaluation
	below is unchanged for a general $M$, with the multiples $\ell N$ replaced by $\ell MN$ and the resulting bound multiplied by
	$M^{-s}$, so the smallest $M$ gives the largest bound. By \eqref{eq:pulse},
	$F_{0}=\sum_{j\in J_{\nu}}\cb{j}\,b\,e^{ijN\pi x}$, and, writing
	$\vartheta=(k+1)\pi/2$,
	\[
	c_{m}=\sum_{j\in J_{\nu}}\cb{j}\int_{0}^{1}b\,\Psi_{m}e^{ijN\pi x},
	\qquad
	\Psi_{m}e^{ijN\pi x}
	=\frac{\sqrt2}{2}\bigl(e^{i\vartheta}e^{i(jN+m)\pi x}
	+e^{-i\vartheta}e^{i(jN-m)\pi x}\bigr),
	\]
	so every integral occurring in $c_{m}$ is a value of the
	transform of $b$,
	\begin{equation}\label{eq:envtr}
		\II(\lambda):=\int_{0}^{1}b(x)e^{i\lambda\pi x}dx
	\end{equation}
	at an exponent $\lambda=jN\pm m$. For $\mu>-1$, $\mu+1$ integrations
	by parts give
	\begin{equation}\label{eq:envtrdecay}
		\abs{\II(\lambda)}\le\frac{\norm{b^{(\mu+1)}}_{L^{1}}}
		{(\pi\abs\lambda)^{\mu+1}}\qquad(\lambda\ne0),
	\end{equation}
	the boundary terms vanishing because $b$ vanishes to order $\mu+1$ at both
	endpoints. Thus $\abs{\II(\lambda)}$ is governed by $\abs\lambda$ alone, and
	we are looking for the leading terms: the exponents $jN\pm m$,
	$j\in J_{\nu}$, are of the order of $N$ unless $m$ lies close to a multiple
	$\ell N$ of $N$, $1\le\ell\le\nu$, so the leading terms of $c_{m}$ come
	from $m$ where the distances $m-\ell N$ are of lower order than $N$. We use
	the width $W$ to separate these distances, which satisfy
	\begin{equation}\label{eq:Wchoice}
		W\to\infty,\qquad W=o(N)\qquad(\text{for instance }W= \lfloor N^{\alpha}\rfloor, 0<\alpha<1 ).
	\end{equation}
	With this width, the coefficients $c_m$ take the following form.

	\begin{lemma}\label{lem:coefs}
		Let $1\le\nu\le k$ and let $W$ satisfy \eqref{eq:Wchoice}. Then,
		uniformly in the cases below,
		\begin{equation}\label{eq:thmcoefs}
			{\small c_{m}=\begin{cases}
				\displaystyle
				\tfrac{\sqrt2}{2}\bigl(e^{-i\vartheta}\cb{\ell}\,\II(-q)
				+e^{i\vartheta}\cb{-\ell}\,\II(q)\bigr)+O\bigl(N^{-\mu-1}\bigr),
				& m=\ell N+q,\ 1\le\ell\le\nu, \ |q|\leq W \\[1mm]
				\displaystyle
				\cb{0}\int_{0}^{1}b\,\Psi_{m}+O\bigl(W^{-\mu-1}\bigr),
				& 1\leq m<N-W,\\[1mm]
				\displaystyle
				O\bigl(W^{-\mu-1}\bigr), & \text{otherwise}.
			\end{cases}}
		\end{equation}
	\end{lemma}

	\begin{proof}
		Recall that $c_{m}$ is a combination of the values $\II(jN\pm m)$ with
		$j\in J_{\nu}$, and that $\II(\lambda)$ is small for large $\abs\lambda$
		by \eqref{eq:envtrdecay}.

		In the first case, $m=\ell N+q$ with $1\le\ell\le\nu$ and $\abs q\le W$.
		If $\pm\ell\in J_{\nu}$ (that is, $\ell$ has the same parity as $\nu$),
		the only small exponents are $jN-m$ at $j=\ell$ and $jN+m$ at $j=-\ell$,
		equal to $-q$ and $q$, respectively, so the two values $\II(\pm q)$ give
		the first case of \eqref{eq:thmcoefs}. If $\pm\ell\notin J_{\nu}$,
		there is no small exponent at all, but the coefficients $\cb{\ell}$ and
		$\cb{-\ell}$ vanish, so the leading term of the first case is zero and
		the same formula holds. In either case the remaining exponents are at
		least $N-W\ge N/2$ in absolute value and contribute $O(N^{-\mu-1})$.
		These are the coefficients that carry the mass near the multiples.

		In the second case, $m<N-W$, every exponent with $j\ne0$ is at least $W$ in
		absolute value, so the terms $j\ne0$ contribute $O(W^{-\mu-1})$ and only the
		term $j=0$ survives, which is the second case of \eqref{eq:thmcoefs}; when
		$\nu$ is odd there is no such term and the estimate holds with $\cb{0}=0$.

		In the third case every exponent is at least $W$ in absolute value, so
		$\abs{c_{m}}\le CW^{-\mu-1}$.
	\end{proof}

	The standard Parseval identities for \eqref{eq:envtr} \cite[Chapter~3]{SS03} give
	\begin{equation*}
		\sum_{q\in\Z}\abs{\II(q)}^{2}=2\norm{b}^{2}_{L^{2}},
		\qquad
		\sum_{q\in\Z}\II(q)^{2}=0,
	\end{equation*}
	hence
	\begin{equation}\label{eq:perio}
		\sum_{q\in\Z}\bigl(\operatorname{Re}\II(q)\bigr)^{2}
		=\sum_{q\in\Z}\bigl(\operatorname{Im}\II(q)\bigr)^{2}=\norm{b}^{2}_{L^{2}}.
	\end{equation}
	We are now ready to give a lower bound for $\lamone^{1/2}$ that is sharper than \eqref{eq:cs} in many cases.

	\begin{lemma}\label{thm:twochan}
		Let $1\le\nu\le k$ and $k\geq s-1$, and let $A_{\nu,r}$ be defined by \eqref{eq:Abump}. Then,
		\begin{equation}\label{eq:twochaneq}
			\lamone^{1/2}\ \ge\ (N\pi)^{-s}\Bigl[
			\frac{\sum_{1\le\abs{j}\le\nu}\abs{\cb{j}}^{2}
				j^{\,2(k+1-s)}}{\Ab{k+1}}
			+\frac{\abs{\cb{0}}^{2}}
			{\Ab{k+1}}\,(N\pi)^{\,2(s-k-1)}
			\Bigr]^{1/2}\bigl(1-o(1)\bigr).
		\end{equation}
		Consequently,
		\begin{equation}\label{eq:lambdonHs}
			\lamone^{1/2}\ \ge\ (N\pi)^{-s} \Bigl(\frac{\Ab{k+1-s}}{\Ab{k+1}}\Bigr)^{1/2}\bigl(1-o(1)\bigr).
		\end{equation}
	\end{lemma}

	\begin{proof}
		By \eqref{eq:onelb} and \eqref{eq:cosnormE},
		\[
		\lamone^{1/2}\ \ge\
		\frac{\norm{F_{0}^{(k+1)}}_{\Hm}}{\norm{F_{0}^{(k+1)}}_{L^{2}}}
		\ \ge\
		\frac{1}{\norm{F_{0}^{(k+1)}}_{L^{2}}}
		\Bigl(\sum_{m\ge1}\frac{(m\pi)^{2(k+1)}}{\Ss(m)}c_{m}^{2}\Bigr)^{1/2},
		\]
		where the $L^{2}$-norm in the denominator is the case $r=k+1$ of
		Lemma~\ref{lem:moments},
		\begin{equation}\label{eq:thmden}
			\norm{F_{0}^{(k+1)}}^{2}_{L^{2}}=(N\pi)^{2(k+1)}\norm{b}^{2}_{L^{2}}\Ab{k+1}
			\bigl(1+O(N^{-2})\bigr).
		\end{equation}
		It remains to estimate the series, and for this we use the three cases of
		\eqref{eq:thmcoefs}; its terms are positive, so the sub-sums of disjoint
		parts add and the third case may be dropped.

		We begin with the first case, $m=\ell N+q$ with $1\le\ell\le\nu$ and
		$\abs q\le W$. For every $\ell$, squaring the first case of \eqref{eq:thmcoefs} gives
		\[
		\sum_{\abs q\le W}c_{\ell N+q}^{2}
		=\frac12\sum_{\abs q\le W}\Bigl|e^{-i\vartheta}\cb{\ell}\II(-q)
		+e^{i\vartheta}\cb{-\ell}\II(q)\Bigr|^{2}+O\bigl(WN^{-\mu-1}\bigr).
		\]
		Since $e^{i\vartheta}=e^{-i\vartheta}(-1)^{k+1}$, $\cb{-\ell}=(-1)^{\nu}\cb{\ell}$, and $\II(-q)=\overline{\II(q)}$, the bracket is
		\[
		\frac14\Bigl|e^{-i\vartheta}\cb{\ell}\II(-q)
		+e^{i\vartheta}\cb{-\ell}\II(q)\Bigr|^{2}=\abs{\cb{\ell}}^{2}\cdot
		\begin{cases}\bigl(\operatorname{Re}\II(q)\bigr)^{2},&\nu+k+1\ \text{even},\\[0.5mm]
			\bigl(\operatorname{Im}\II(q)\bigr)^{2},&\nu+k+1\ \text{odd},\end{cases}
		\]
		By \eqref{eq:perio}, the two alternatives give the same square-sum
		$\norm{b}^{2}_{L^{2}}$ with the truncation error
		\begin{equation*}
			\sum_{\abs q>W}\abs{\II(q)}^{2}\le C\sum_{\abs q>W}\abs q^{-2\mu-2}
			=O\bigl(W^{-2\mu-1}\bigr).
		\end{equation*}
		Hence, for every $\ell$,
		\begin{equation}\label{eq:blockmass}
			\sum_{\abs q\le W}c_{\ell N+q}^{2}
			=2\abs{\cb{\ell}}^{2}\norm{b}^{2}_{L^{2}}
			\bigl(1+O(W^{-2\mu-1}+WN^{-\mu-1})\bigr).
		\end{equation}
		Since $\Ss(m)=(m\pi)^{2s}(1+O(m^{-2}))$ and $m=\ell N+q$ with
		$\abs q\le W=o(N)$,
		\[
		\frac{(m\pi)^{2(k+1)}}{\Ss(m)}
		=(m\pi)^{2(k+1-s)}\bigl(1+O(m^{-2})\bigr)
		=(\ell N\pi)^{2(k+1-s)}\bigl(1+O(WN^{-1})\bigr).
		\]
		Combining \eqref{eq:blockmass} and \eqref{eq:thmden}, we obtain
		\begin{align*}
			& \frac{1}{\norm{F_{0}^{(k+1)}}^2_{L^{2}}}
			\sum_{m=\ell N+q,\ \abs q\le W}\frac{(m\pi)^{2(k+1)}}{\Ss(m)}c_{m}^{2} \\
			&= 2\frac{\abs{\cb{\ell}}^2\ell^{\,2(k+1-s)}(N\pi)^{-2s}}{\Ab{k+1}} + O(W^{-2\mu-1}+WN^{-\mu-1}+WN^{-1}).
		\end{align*}
		Summing over the multiples, that is, over the values $1\le\abs j\le\nu$
		(the two signs $j=\pm\ell$ contribute to the same block), we obtain
		\begin{align*}
			&\frac{1}{\norm{F_{0}^{(k+1)}}^{2}_{L^{2}}}
			\sum_{\substack{m=\ell N+q\\ 1\le \ell\le\nu,|q|\le W }} \frac{(m\pi)^{2(k+1)}}{\Ss(m)}c_{m}^{2} \nonumber\\
			&\ge\ (N\pi)^{-2s}
			\frac{\sum_{1\le\abs j\le\nu}\abs{\cb{j}}^{2}j^{\,2(k+1-s)}}{\Ab{k+1}}-O(W^{-2\mu-1}+WN^{-\mu-1}+WN^{-1}).
		\end{align*}
		This is the first term of \eqref{eq:twochaneq}.

		We compute the second case, $1\le m<N-W$, next. By the second case of
		\eqref{eq:thmcoefs}, $c_{m}=\cb{0}\int_{0}^{1}b\,\Psi_{m}+r_{m}$ with
		$r_{m}= O(W^{-\mu-1})$. Squaring gives
		\[
		c_{m}^{2}
		=\abs{\cb{0}}^{2}\Bigl(\int_{0}^{1}b\,\Psi_{m}\Bigr)^{2}
		+2\cb{0}\Bigl(\int_{0}^{1}b\,\Psi_{m}\Bigr)r_{m}
		+r_{m}^{2}.
		\]
		We sum the first two terms against the weight $(m\pi)^{2(k+1)}/\Ss(m)$
		over the values $1\le m<N-W$. The third  term is positive and can be dropped. The first
		contribution is
		\begin{equation}\label{eq:low-mode}
			\abs{\cb{0}}^{2}\sum_{m<N-W}\frac{(m\pi)^{2(k+1)}}{S_s(m)}
			\Bigl(\int_{0}^{1}b\,\Psi_{m}\Bigr)^{2}.
		\end{equation}
		The second contribution is at most
		\begin{align*}
			&2\abs{\cb{0}}\Bigl(\sup_{m<N-W}\abs{r_{m}}\Bigr)\Bigl(\max_{m<N-W}
			\frac{(m\pi)^{2(k+1)}}{\Ss(m)}\Bigr)
			\sum_{m<N-W}\Bigl|\int_{0}^{1}b\,\Psi_{m}\Bigr| \\
			&\le\ 2\abs{\cb{0}}\Bigl(\sup_{m<N-W}\abs{r_{m}}\Bigr)(N\pi)^{2(k+1-s)} \sum_{m<N-W}\Bigl|\int_{0}^{1}b\,\Psi_{m}\Bigr|.
		\end{align*}
		Since $b^{(l)}(0)=b^{(l)}(1)=0$ for $0\le l\le\mu$, $\{\Psi_m\}_{m\ge1}=\pm\sqrt2\{\sin m\pi x\}_{m\ge1}$ for $k$ even, and
		$\{\Psi_m\}_{m\ge1}=\pm\sqrt2\{\cos m\pi x\}_{m\ge1}$ for $k$ odd, integration by parts $\mu+1$ times gives
		\[
		\Bigl|\int_{0}^{1}b\,\Psi_{m}\Bigr|=O\bigl(m^{-(\mu+1)}\bigr).
		\]
		Hence
		\[
		\sum_{m<N-W}\Bigl|\int_{0}^{1}b\,\Psi_{m}\Bigr|
		=\begin{cases}
			O(\log N),&\mu=0,\\
			O(1),&\mu\ge1.
		\end{cases}
		\]
		After division by \eqref{eq:thmden}, the second contribution is therefore
		$O\bigl(\sup_{m<N-W}\abs{r_m}\,N^{-2s}\log N\bigr)=o(N^{-2s})$, because $\sup_{m<N-W}\abs{r_m}=O(W^{-\mu-1})$
		and $W\gg\log N$ by taking $W=\lfloor N^{\alpha}\rfloor$, $0<\alpha<1$. Thus it does not affect the leading term $N^{-2s}$.

		What remains is to analyze \eqref{eq:low-mode}, whose weight is
		\begin{equation}\label{eq:lowweight}
			\frac{(m\pi)^{2(k+1)}}{\Ss(m)}
			=\frac{(m\pi)^{2(k+1-s)}}{1+\delta_{m}},
			\qquad
			\delta_{m}:=\sum_{a=0}^{s-1}(m\pi)^{2(a-s)}\ge0.
		\end{equation}
		We can take $b$ to be sufficiently oscillatory so that the integrals
		$\int_{0}^{1}b\,\Psi_{m}$ are small away from $m=K$ with $1\ll K\ll W$ (e.g. $K=\lfloor \log N \rfloor$, $W=\lfloor N^{\alpha}\rfloor$, $0<\alpha<1$). Indeed, let
		\[
		b:=b_{0}\bigl(\sin(K\pi x)-\lambda\bigr),\qquad
		b_{0}:=x^{\mu+1}(1-x)^{\mu+1},
		\]
		with $\lambda$ chosen so that $\int_{0}^{1}b=0$. Assume first that $k$ is
		even (the case of odd $k$ is analogous). Then
		$\Psi_{m}=\pm\sqrt2\sin m\pi x$, and
		\[
		\int_{0}^{1}b\,\Psi_{m}
		=\pm\frac{\sqrt2}{2}\int_{0}^{1}b_{0}\bigl(\cos((K-m)\pi x)
		-\cos((K+m)\pi x)\bigr)+O\bigl(K^{-\mu-1}m^{-\mu-1}\bigr),
		\]
		where the error term is the contribution
		$-\lambda\int_{0}^{1}b_{0}\Psi_{m}$, with
		$\lambda=\int_{0}^{1}b_{0}\sin K\pi x/\int_{0}^{1}b_{0}
		=O(K^{-\mu-1})$. The integrals are large only near $m=K$. At $m=K$,
		\begin{equation*}
			\int_{0}^{1}b\,\Psi_{K} = \pm\frac{\sqrt2}{2}\int_{0}^{1}b_{0} + O(K^{-\mu-1})
		\end{equation*}
		is not small. Away from $m=K$, we have the smallness
		\[
		\sum_{\abs{m-K}>K/2}\Bigl(\int_{0}^{1}b\,\Psi_{m}\Bigr)^{2}\leq C\sum_{\abs{m-K}>K/2}\abs{m-K}^{-2\mu-2}
		\ =\ O\bigl(K^{-2\mu-1}\bigr).
		\]
		With this choice of $b$, we show that \eqref{eq:low-mode} is asymptotic to the series
		\[
		\abs{\cb{0}}^{2}\sum_{m<N-W}(m\pi)^{2(k+1-s)}
		\Bigl(\int_{0}^{1}b\,\Psi_{m}\Bigr)^{2}:= \abs{\cb{0}}^{2} S .
		\]
		We compare the two on the two ranges of $m$. For $K/2<m<N-W$, by \eqref{eq:lowweight} and $\delta_{m}=O(K^{-2})$ in this range, we have
		\begin{equation}\label{eq:weight-bound}
			\sum_{K/2<m<N-W}\frac{(m\pi)^{2(k+1)}}{\Ss(m)}
			\Bigl(\int_{0}^{1}b\,\Psi_{m}\Bigr)^{2}
			=\sum_{K/2<m<N-W}(m\pi)^{2(k+1-s)}
			\bigl(1+O(K^{-2})\bigr)\Bigl(\int_{0}^{1}b\,\Psi_{m}\Bigr)^{2}.
		\end{equation}
		For $m\le K/2$, $\delta_{m}$ is only $O(1)$, while $\int_{0}^{1}b\,\Psi_{m}$ is small. Notice that
		\[
		S_1:=\sum_{m\le K/2}(m\pi)^{2(k+1-s)}
		\Bigl(\int_{0}^{1}b\,\Psi_{m}\Bigr)^{2} \leq CK^{2(k+1-s)}\sum_{m\le K/2}\Bigl(\int_{0}^{1}b\,\Psi_{m}\Bigr)^{2}\leq CK^{2(k+1-s)} K^{-2\mu-1},
		\]
		and
		\begin{align*}
			S_2:=\sum_{K/2< m< N-W}(m\pi)^{2(k+1-s)}
			\Bigl(\int_{0}^{1}b\,\Psi_{m}\Bigr)^{2} &\geq CK^{2(k+1-s)}\sum_{K/2< m< N-W}\Bigl(\int_{0}^{1}b\,\Psi_{m}\Bigr)^{2}\\
			&\geq CK^{2(k+1-s)}\Bigl(\int_{0}^{1}b\,\Psi_{K}\Bigr)^{2}.
		\end{align*}
		Then, $(1+O(K^{-2\mu-1}))S_2 = S_1+S_2=S$, and by \eqref{eq:weight-bound},
		\[
		\sum_{K/2<m<N-W}\frac{(m\pi)^{2(k+1)}}{\Ss(m)}
		\Bigl(\int_{0}^{1}b\,\Psi_{m}\Bigr)^{2} = (1+O(K^{-2}))S_2 = (1+O(K^{-2})+O(K^{-2\mu-1}))S.
		\]
		By \eqref{eq:lowweight} with $\delta_m\geq0$,
		\begin{equation*}
			\sum_{m\le K/2}\frac{(m\pi)^{2(k+1)}}{\Ss(m)}
			\Bigl(\int_{0}^{1}b\,\Psi_{m}\Bigr)^{2}
			\ \le\ \sum_{m\le K/2}(m\pi)^{2(k+1-s)}
			\Bigl(\int_{0}^{1}b\,\Psi_{m}\Bigr)^{2} = S_1 = O(K^{-2\mu-1})S.
		\end{equation*}
		Hence \eqref{eq:low-mode} is at least
		\begin{equation*}
			\abs{\cb{0}}^{2}(1-O(K^{-2})-O(K^{-2\mu-1}))S.
		\end{equation*}
		By the tail estimate ($\abs{\int_{0}^{1}b\Psi_{m}}
		=O(m^{-\mu-1})$, the sum $\sum_{m<N-W}
		\Bigl(\int_{0}^{1}b\,\Psi_{m}\Bigr)^{2} $ differs
		from the full sum by $O(N^{-2\mu-1})$); hence
		\begin{align*}
			\eqref{eq:low-mode}\ \ge&\ \abs{\cb{0}}^{2}(1-O(K^{-2})-O(K^{-2\mu-1}))\sum_{m<N-W}
			\Bigl(\int_{0}^{1}b\,\Psi_{m}\Bigr)^{2} \\
			\ge&\ \ \abs{\cb{0}}^{2}(1-O(K^{-2})-O(K^{-2\mu-1})- O(N^{-2\mu-1}))\sum_{m\geq 1}
			\Bigl(\int_{0}^{1}b\,\Psi_{m}\Bigr)^{2}\\
			\ge&\
			\abs{\cb{0}}^{2}\norm{b}^{2}_{L^{2}}\bigl(1-o(1)\bigr),
		\end{align*}
		the $o(1)$ being $O(K^{-2})+O(K^{-2\mu-1})+O(N^{-2\mu-1})$. Here, the last inequality uses the fact
		\begin{equation*}
			\sum_{m\ge1}\Bigl(\int_{0}^{1}b\,\Psi_{m}\Bigr)^{2}
			=\begin{cases}
				\sum_{m\ge1}\bigl(\int_{0}^{1}b\,\sqrt2\sin m\pi x\bigr)^{2}
				=\norm{b}^{2}_{L^{2}},&k\ \text{even},\\[0.6mm]
				\sum_{m\ge1}\bigl(\int_{0}^{1}b\,\sqrt2\cos m\pi x\bigr)^{2}
				=\norm{b}^{2}_{L^{2}}-\bigl(\int_{0}^{1}b\bigr)^{2},&k\ \text{odd}.
			\end{cases}
		\end{equation*}
		Dividing by \eqref{eq:thmden} and using the splitting
		$(N\pi)^{-2(k+1)}=(N\pi)^{2(s-k-1)}(N\pi)^{-2s}$, \eqref{eq:low-mode} becomes the second
		term of \eqref{eq:twochaneq}.
	\end{proof}

\begin{remark}\label{rem:nuplusone}
	Lemma~\ref{lem:coefs} and Lemma~\ref{thm:twochan} are stated for
	$1\le\nu\le k$ because the estimate \eqref{eq:envtrdecay}, which depends on
	the boundary condition imposed on $b$, decays in $\abs\lambda$ only for
	$\mu>-1$. It suffices to assume that $b$ vanishes to one order more at the
	endpoints, $b^{(l)}(0)=b^{(l)}(1)=0$ for $l\le\mu+1$ (for $\nu=k+1$ this
	just says $b(0)=b(1)=0$), and such a choice is admissible. Then
	\eqref{eq:envtrdecay} is governed by $\abs\lambda^{-(\mu+2)}$, which tends
	to $0$ for $\mu=-1$. The statements of Lemma~\ref{lem:coefs} and of
	Lemma~\ref{thm:twochan} are therefore valid for $\nu=k+1$ as well.
\end{remark}

\begin{remark}\label{rem:constantK}
	In the proof of Lemma~\ref{thm:twochan}, $b$ is taken oscillatory, which
	makes the constant $\norm{b^{(\mu+1)}}_{L^{1}}$ of \eqref{eq:envtrdecay}
	large, of order $K^{\mu+1}$, and no longer a constant. However, since
	$K\ll W\ll N$, all error terms remain governed by $W$ and $N$. The choice
	of $b$ is therefore legitimate.
\end{remark}

	Dividing the lower bound \eqref{eq:lambdonHs} of Lemma~\ref{thm:twochan} by
	the $n$-width \eqref{eq:dnform} gives the following lower bound for
	the ratio $\Lam$.

	\begin{theorem}\label{thm:nusharp}
			Let $1\le\nu\le k+1$ and $k\ge s-1$. Then,
		\begin{equation}\label{eq:nusharp}
			\Lam\ \ge\ \nu^{s}\Bigl(\frac{\Ab{k+1-s}}{\Ab{k+1}}\Bigr)^{1/2}
			\bigl(1-o(1)\bigr),
		\end{equation}
		where $A_{\nu,r}$ is defined by \eqref{eq:Abump}.
	\end{theorem}

The bound \eqref{eq:nusharp} is stronger than \eqref{eq:nucrude} of the previous subsection. In fact, splitting $|c_{\nu,j}|^2j^{2k+2}=|c_{\nu,j}| j^{k+1-s}|c_{\nu,j}| j^{k+1+s}$ and applying Cauchy--Schwarz to $A_{\nu,k+1}$ gives
\[
	A_{\nu,k+1}\le\Bigl(\sum_{j\in J_{\nu}}|c_{\nu,j}| ^{2} j^{2(k+1-s)}\Bigr)^{1/2}
\Bigl(\sum_{j\in J_{\nu}}|c_{\nu,j}| ^{2} j^{2(k+1+s)}\Bigr)^{1/2}
= A_{\nu,k+1-s}^{1/2}A_{\nu,k+1+s}^{1/2}.
\]
Equality in Cauchy--Schwarz holds exactly when the nonzero coefficients
$c_{\nu,j}$ are supported on a single value of $\abs j$. Here this happens for
$\nu=1$ and, when $s\le k$, for $\nu=2$; in these cases, and only in these, the two lower bounds coincide.
At the borderline $k=s-1$, the new bound \eqref{eq:nusharp} is strictly larger than $1$ for every $\nu\ge2$. For $\nu=2$ this is new information, since the bound \eqref{eq:nucrude} equals $1$ there. Indeed, by \eqref{eq:Abump},
\[
	\nu^{2s}\Bigl(\frac{\Ab{0}}{\Ab{s}}\Bigr)-1
	=\frac{\sum_{m=0}^{\nu}\binom{\nu}{m}^{2}\bigl(\nu^{2s}-\abs{\nu-2m}^{2s}\bigr)}
	{\sum_{m=0}^{\nu}\binom{\nu}{m}^{2}\abs{\nu-2m}^{2s}}>0,
\]
every numerator term being non-negative and the term $m=1$ positive for $\nu\ge2$. For the most regular ball, with $s=k+1$, the maximally smooth spline is thus strictly superior, in the sense of $n$-widths and at equal dimension, to every less smooth spline ($\nu\ge2$).

By contrast, for $k>s-1$ and $\mu=k-2$ ($\nu=2$) both bounds reduce to the value $1$: with $A_{2,r}=2^{2r-3}$ and $k+1-s\ge1$,
\[
	\nu^{s}\Bigl(\frac{\Ab{k+1-s}}{\Ab{k+1}}\Bigr)^{1/2}=2^{s}\bigl(2^{-2s}\bigr)^{1/2}=1,
\]
just as in \eqref{eq:nucrude}. The numerical experiments of Section~\ref{sec:numerics} indicate that the value $1$ is sharp here, that is,  the ratio $\Lam$ tends to $1$ for this space. For targets that are not maximally smooth ($s<k+1$), the approximation in the $B$-spline space $\Kh$ and in the space with $\mu=k-2$ is asymptotically the same, and both reach the $n$-width at equal dimension. To the authors' knowledge, this asymptotic coincidence has not been recorded in the earlier literature.

We close this subsection with the following corollary.

\begin{corollary}\label{cor:dichotomy}
	Let $1\le\nu\le k+1$ and $k\ge s-1$. Then, for a sufficiently large $n$,
	\begin{itemize}\setlength{\itemsep}{0pt}
		\item If $k=s-1$, then $\Lam>1$ for every $\nu\ge2$;
		\item If $k>s-1$, then $\Lam>1$ for every $\nu>2$, while for $\nu=2$, the constant factor in the lower bound equals $1$.
	\end{itemize}
\end{corollary}

	\section{Numerical experiments}\label{sec:numerics}

	In this section, we numerically approximate the $n$-width $\dn$ and the
	sup-inf $E(\Bball,\Xsp;L^{2})$, and through them the ratio $\Lam$. Both
	quantities are read from eigenvalue problems, whose eigenvalues we
	approximate by the Galerkin method. The two variational problems are
	discretized in a finite-dimensional subspace, in the form used in
	\cite[Subsection 5.1]{EBH09}.

	The two problems are discretized in one and the same subspace. Choose
	$X_{m}\subset\Hsp$ with $\dim X_{m}=m$ and $m\gg n$. The width $\dn$ is
	the square root of the $(n+1)$-st eigenvalue of the free problem
	\eqref{eq:free}, in the ordering of Lemma~\ref{lem:dneig}. Its Galerkin
	discretization is the finite-dimensional eigenproblem of finding
	$u\in X_{m}$ and $z\in\mathbb{R}$ such that
	\begin{equation}\label{eq:xmfree}
		(u,v)_{L^{2}}=z\,(u,v)_{\Hsp}\qquad \forall v\in X_{m}.
	\end{equation}
	Likewise, $E(\Bball,\Xsp;L^{2})$ is the square root of the largest
	eigenvalue of \eqref{eq:supeig}. Its Galerkin discretization is the problem
	of finding $\tilde u\in X_{m}$ and $\tilde z\in\mathbb{R}$ such that
	\begin{equation}\label{eq:xmE}
		\bigl((I-\Pi_{n})\tilde u,\,(I-\Pi_{n})v\bigr)_{L^{2}}
		=\tilde z\,(\tilde u,v)_{\Hsp}\qquad \forall v\in X_{m}.
	\end{equation}

	We rewrite \eqref{eq:xmfree} and \eqref{eq:xmE} as matrix eigenproblems.
	Let $\{N_{0,i}\}_{i=1}^{n}$ be a basis of $\Xsp$ and $\{N_{s,i}\}_{i=1}^{m}$
	a basis of $X_{m}$, and for $u=\sum_{i}u_{i}N_{s,i}$ write
	$\boldsymbol u=(u_{i})$ for its coefficient vector. Define the matrices
	\begin{equation*}
		\begin{aligned}
			H_{i,j}&=(N_{0,i},N_{0,j})_{L^{2}},\qquad
			M_{i,j}=(N_{s,i},N_{s,j})_{L^{2}},\\
			K_{i,j}&=(N_{s,i},N_{s,j})_{\Hsp},\qquad
			T_{i,j}=(N_{0,i},N_{s,j})_{L^{2}},
		\end{aligned}
	\end{equation*}
	of sizes $n\times n$, $m\times m$, $m\times m$ and $n\times m$, and the
	Schur complement
	\begin{equation*}
		\tilde M:=M-T^{\mathsf T}H^{-1}T .
	\end{equation*}
	The matrices $M$ and $K$ are the mass and the energy matrix of $X_{m}$,
	hence positive definite, and $H$ is the mass matrix of $\Xsp$; the matrix
	$\tilde M$ is the matrix of the form
	$\bigl((I-\Pi_{n})u,\,(I-\Pi_{n})v\bigr)_{L^{2}}$ restricted to $X_{m}$.
	Indeed, the coefficients of the projection $\Pi_{n}u=\sum_{i}c_{i}N_{0,i}$
	are characterized by $(u-\Pi_{n}u,N_{0,i})_{L^{2}}=0$ for all $i$, that is,
	by the normal equations
	\begin{equation*}
		H\boldsymbol c=T\boldsymbol u,\qquad
		(T\boldsymbol u)_{i}=(u,N_{0,i})_{L^{2}}=\sum_{j}T_{i,j}u_{j},
	\end{equation*}
	so that $(\Pi_{n}u,\Pi_{n}v)_{L^{2}}
	=\boldsymbol c(u)^{\mathsf T}H\boldsymbol c(v)
	=(T\boldsymbol u)^{\mathsf T}H^{-1}T\boldsymbol v$; subtracting
	$(u,v)_{L^{2}}=\boldsymbol u^{\mathsf T}M\boldsymbol v$, the matrix of the
	form is $M-T^{\mathsf T}H^{-1}T=\tilde M$, and in particular
	$\boldsymbol u^{\mathsf T}\tilde M\boldsymbol u
	=\norm{(I-\Pi_{n})u}^{2}_{L^{2}}\ge0$, with equality exactly when
	$u\in\Xsp$. Thus $\tilde M$ is positive semidefinite, with null space
	$X_{m}\cap\Xsp$. Let $\tilde n = \dim (X_{m}\cap\Xsp)$. The
	two eigenproblems \eqref{eq:xmfree} and \eqref{eq:xmE} are equivalent to
	\begin{equation}\label{eq:83}
		Mu=zKu
	\end{equation}
	and
	\begin{equation}\label{eq:84}
		\tilde Mu=\tilde zKu,
	\end{equation}
	respectively.
Their eigenvalues are the approximations of the two quantities,
	\begin{equation*}
		\dn\approx z_{n+1}^{1/2},\qquad
		E(\Bball,\Xsp;L^{2})\approx\tilde z_{1}^{1/2}
	\end{equation*}
	where $z_{1}\ge z_{2}\ge\cdots\ge z_{m}>0$ are the eigenvalues of
	\eqref{eq:83}, and $\tilde z_{1}\ge\cdots\ge\tilde z_{m-\tilde n}>0$ are
	the positive eigenvalues of \eqref{eq:84}.
	Since $X_{m}\subset\Hsp$, the Rayleigh--Ritz principle gives
	$\tilde z_{1}\le\tilde\lambda_{1}$ and $z_{j}\le\lambda_{j}$ for every
	$j$. The Galerkin methods converge as $m\to\infty$.

	In the computations, $X_{m}$ is taken as follows. Its mesh is a nested
	refinement of the coarse mesh of $\Xsp$: each of the $N$ coarse elements
	is divided into an integer number of equal parts, so that $m$ is increased with the number of parts. Let $X_{m}$ be the space of maximally smooth
	splines of degree $k^{*}$ on those meshes, and set
	$k^{*}=\max\{s+3,7\}$, so as to improve the rate of approximation as the
	mesh is refined. For each eigenvalue problem, the size $m$ of $X_{m}$ was
	increased by refining the nested mesh, until the computed values were
	stable in the digits shown.
	The numerical results are summarized as follows.
	\begin{table}[htbp]
		\centering
		\setlength{\tabcolsep}{5pt}
		\caption{Numerical approximation of the $n$-width through the eigenvalues $z_{n+1}$, with $q_{n}:=z_{n+1}^{-1/(2s)}/\pi$. The table is consistent with $\abs{n-q_{n}}\to(s-1)/2$ and with boundedness of $\bigl|n\bigl(n-(s-1)/2-q_{n}\bigr)\bigr|$.}
		\label{tab:1}
		\begin{tabular}{c|ccc|ccc}
			\multirow{2}{*}{$s$} & \multicolumn{3}{c|}{$|n-q_{n}|$} & \multicolumn{3}{c}{$\bigl|n\bigl(n-(s-1)/2-q_{n}\bigr)\bigr|$}\\
			& $n=20$ & $n=40$ & $n=60$ & $n=20$ & $n=40$ & $n=60$\\ \hline
			$1$ & $0.0025329$ & $0.0012665$ & $0.00084434$ & $0.050657$ & $0.050660$ & $0.050660$\\
			$2$ & $0.49857$ & $0.49933$ & $0.49956$ & $0.028528$ & $0.026892$ & $0.026363$\\
			$3$ & $0.99885$ & $0.99951$ & $0.99969$ & $0.022936$ & $0.019769$ & $0.018778$\\
			$4$ & $1.4989$ & $1.4996$ & $1.4997$ & $0.021656$ & $0.016836$ & $0.015379$\\
			$5$ & $1.9989$ & $1.9996$ & $1.9998$ & $0.022287$ & $0.015615$ & $0.013668$\\
			$6$ & $2.4988$ & $2.4996$ & $2.4998$ & $0.024060$ & $0.015283$ & $0.012815$\\
		\end{tabular}
	\end{table}
	\begin{table}[ht]
		\centering
		\setlength{\tabcolsep}{4pt}
		\caption{Numerical approximation of the sup-inf error from $\tilde z_{1}$ for some pairs $(\nu,k,s)$ with $\nu\le k$. Here $r_{1}=(\Ab{k+1}/\Ab{k+1+s})^{1/2}$, $r_{2}=(\Ab{k+1-s}/\Ab{k+1})^{1/2}$, $D$ is the ratio $r_{1}/r_{2}$, and $\theta_{N}=(N\pi)^{s}\tilde z_{1}^{1/2}$.}
		\label{tab:2}
		\begin{tabular}{c|c|c|c|ccc|c}
			\multirow{2}{*}{$(\nu,k,s)$} & \multirow{2}{*}{$r_1$} & \multirow{2}{*}{$r_2$} & \multirow{2}{*}{$D$} & \multicolumn{3}{c|}{$\theta_{N}$} & \multirow{2}{*}{$\theta_{60}/r_2$}\\
			& & & & $N=20$ & $N=40$ & $N=60$ & \\ \hline
			(1,1,1) & 1.0000 & 1.0000 & 1.0000 & 0.99377 & 0.99843 & 0.99930 & 0.9993\\
			(1,1,2) & 1.0000 & 1.0000 & 1.0000 & 0.99275 & 0.99817 & 0.99919 & 0.9992\\
			(1,2,2) & 1.0000 & 1.0000 & 1.0000 & 0.98172 & 0.99549 & 0.99801 & 0.9980\\
			(1,2,3) & 1.0000 & 1.0000 & 1.0000 & 0.98048 & 0.99518 & 0.99787 & 0.9979\\
			(1,3,3) & 1.0000 & 1.0000 & 1.0000 & 0.96300 & 0.99104 & 0.99608 & 0.9961\\
			(1,3,4) & 1.0000 & 1.0000 & 1.0000 & 0.96170 & 0.99072 & 0.99594 & 0.9959\\
			(1,4,4) & 1.0000 & 1.0000 & 1.0000 & 0.93656 & 0.98491 & 0.99347 & 0.9935\\
			(1,4,5) & 1.0000 & 1.0000 & 1.0000 & 0.93521 & 0.98458 & 0.99332 & 0.9933\\
			(1,5,5) & 1.0000 & 1.0000 & 1.0000 & 0.90139 & 0.97689 & 0.99010 & 0.9901\\
			(1,5,6) & 1.0000 & 1.0000 & 1.0000 & 0.90000 & 0.97655 & 0.98995 & 0.9900\\
			(2,2,1) & 0.50000 & 0.50000 & 1.0000 & 0.49865 & 0.49966 & 0.49985 & 0.9997\\
			(3,3,1) & 0.33354 & 0.33515 & 0.9952 & 0.33466 & 0.33503 & 0.33510 & 0.9998\\
			(2,2,2) & 0.25000 & 0.25000 & 1.0000 & 0.24917 & 0.24979 & 0.24991 & 0.9996\\
			(2,3,2) & 0.25000 & 0.25000 & 1.0000 & 0.24823 & 0.24956 & 0.24981 & 0.9992\\
			(2,3,3) & 0.12500 & 0.12500 & 1.0000 & 0.12401 & 0.12476 & 0.12489 & 0.9991\\
			(2,2,3) & 0.12500 & 0.21651 & 0.5774 & 0.21211 & 0.21463 & 0.21566 & 0.9961\\
			(2,4,4) & 0.062500 & 0.062500 & 1.0000 & 0.061616 & 0.062286 & 0.062406 & 0.9985\\
			(3,3,4) & 0.012354 & 0.039014 & 0.3167 & 0.039757 & 0.039846 & 0.039862 & 1.0218\\
			(3,4,5) & 0.0041155 & 0.013012 & 0.3163 & 0.013028 & 0.013088 & 0.013100 & 1.0067\\
			(4,5,5) & 0.00097847 & 0.0021794 & 0.4490 & 0.0022007 & 0.0022076 & 0.0022090 & 1.0136\\
			(2,3,4) & 0.062500 & 0.10825 & 0.5774 & 0.10399 & 0.10651 & 0.10702 & 0.9886\\
			(4,4,5) & 0.00098415 & 0.0057328 & 0.1717 & 0.0059874 & 0.0060335 & 0.0060506 & 1.0554\\
			(2,5,6) & 0.015625 & 0.027063 & 0.5774 & 0.025604 & 0.026259 & 0.026582 & 0.9822\\
			(3,5,6) & 0.0013718 & 0.0043378 & 0.3162 & 0.0043068 & 0.0043386 & 0.0043443 & 1.0015\\
			(4,5,6) & 0.00024462 & 0.0014415 & 0.1697 & 0.0014410 & 0.0014609 & 0.0014663 & 1.0172\\
			(5,5,6) & 0.000065714 & 0.00069961 & 0.0939 & 0.00078484 & 0.00078601 & 0.00078623 & 1.1238\\
		\end{tabular}
	\end{table}

	\begin{table}[ht]
		\centering
		\setlength{\tabcolsep}{4pt}
		\caption{Numerical approximation of the sup-inf error from $\tilde z_{1}$ for some pairs $(\nu,k,s)$ with $\nu=k+1$. Here $r_{1}=(\Ab{k+1}/\Ab{k+1+s})^{1/2}$, $r_{2}=(\Ab{k+1-s}/\Ab{k+1})^{1/2}$, $D$ is the ratio $r_{1}/r_{2}$, and $\theta_{N}=(N\pi)^{s}\tilde z_{1}^{1/2}$.}
		\label{tab:3}
		\begin{tabular}{c|c|c|c|ccc|c}
			\multirow{2}{*}{$(\nu,k,s)$} & \multirow{2}{*}{$r_1$} & \multirow{2}{*}{$r_2$} & \multirow{2}{*}{$D$} & \multicolumn{3}{c|}{$\theta_{N}$} & \multirow{2}{*}{$\theta_{60}/r_2$}\\
			& & & & $N=20$ & $N=40$ & $N=60$ & \\ \hline
			(1,0,1) & 1.0000 & 1.0000 & 1.0000 & 0.99987 & 0.99997 & 0.99999 & 1.0000\\
			(2,1,1) & 0.50000 & 0.50000 & 1.0000 & 0.49998 & 0.50000 & 0.50000 & 1.0000\\
			(2,1,2) & 0.25000 & 0.43301 & 0.5774 & 0.44106 & 0.44111 & 0.44112 & 1.0187\\
			(3,2,2) & 0.11179 & 0.15617 & 0.7158 & 0.16002 & 0.16003 & 0.16003 & 1.0247\\
			(3,2,3) & 0.037265 & 0.11641 & 0.3201 & 0.12499 & 0.12500 & 0.12500 & 1.0738\\
			(4,3,3) & 0.016098 & 0.033895 & 0.4749 & 0.037009 & 0.037010 & 0.037010 & 1.0919\\
			(4,3,4) & 0.0040260 & 0.022420 & 0.1796 & 0.026062 & 0.026064 & 0.026065 & 1.1626\\
			(5,4,4) & 0.0017145 & 0.0055797 & 0.3073 & 0.0066631 & 0.0066632 & 0.0066632 & 1.1942\\
			(5,4,5) & 0.00034318 & 0.0033478 & 0.1025 & 0.0042974 & 0.0042977 & 0.0042978 & 1.2837\\
			(6,5,5) & 0.00014503 & 0.00073726 & 0.1967 & 0.00097986 & 0.00097988 & 0.00097988 & 1.3291\\
			(6,5,6) & 0.000024203 & 0.00040754 & 0.0594 & 0.00058586 & 0.00058591 & 0.00058592 & 1.4377\\
			(7,6,6) & 0.000010175 & 0.000081265 & 0.1252 & 0.00012166 & 0.00012166 & 0.00012166 & 1.4971\\
		\end{tabular}
	\end{table}

	\begin{table}[ht]
		\centering
		\setlength{\tabcolsep}{4pt}
		\caption{Numerical approximation of the sup-inf error from $\tilde z_{1}$ for some pairs $(\nu,k,s)$ with $k=s=6$. Here $r_{1}=(\Ab{k+1}/\Ab{k+1+s})^{1/2}$, $r_{2}=(\Ab{k+1-s}/\Ab{k+1})^{1/2}$, $D$ is the ratio $r_{1}/r_{2}$, and $\theta_{N}=(N\pi)^{s}\tilde z_{1}^{1/2}$.}
		\label{tab:4}
		\begin{tabular}{c|c|c|c|ccc|c}
			\multirow{2}{*}{$(\nu,k,s)$} & \multirow{2}{*}{$r_1$} & \multirow{2}{*}{$r_2$} & \multirow{2}{*}{$D$} & \multicolumn{3}{c|}{$\theta_{N}$} & \multirow{2}{*}{$\theta_{60}/r_2$}\\
			& & & & $N=20$ & $N=40$ & $N=60$ & \\ \hline
			(1,6,6) & 1.0000 & 1.0000 & 1.0000 & 0.85670 & 0.96676 & 0.98591 & 0.9859\\
			(2,6,6) & 0.015625 & 0.015625 & 1.0000 & 0.015122 & 0.015508 & 0.015574 & 0.9968\\
			(3,6,6) & 0.0013717 & 0.0019399 & 0.7071 & 0.0019149 & 0.0019344 & 0.0019378 & 0.9989\\
			(4,6,6) & 0.00024426 & 0.00054565 & 0.4477 & 0.00054523 & 0.00054799 & 0.00054848 & 1.0052\\
			(5,6,6) & 0.000064622 & 0.00023715 & 0.2725 & 0.00024591 & 0.00024655 & 0.00024667 & 1.0401\\
			(6,6,6) & 0.000022706 & 0.00013106 & 0.1733 & 0.00015241 & 0.00015258 & 0.00015261 & 1.1645\\
		\end{tabular}
	\end{table}

		\begin{itemize}
			\item Table~\ref{tab:1} illustrates Conjecture~\ref{conj:shift}. For
			$s=1,2$ the statement is covered by Lemma~\ref{lem:spectrum}, so the
			table illustrates the lemma numerically; for $s=3,\dots,6$ it
			supports the conjecture. For $s=1,\dots,6$ and $n=20,40,60$, the
			table lists, for $q_{n}:=z_{n+1}^{-1/(2s)}/\pi$,
			\begin{equation*}
				|n-q_{n}|,\qquad \bigl|n\Bigl(n-\frac{s-1}{2}-q_{n}\Bigr)\bigr|.
			\end{equation*}
			The first block approaches $(s-1)/2$,  with the deviations from the limit decreasing monotonically. The second block stays bounded, which is
			consistent with the $O(n^{-1})$ tail estimate of
			Conjecture~\ref{conj:shift}.
			\item Tables~\ref{tab:2}--\ref{tab:4} report the normalized error
				$\theta_{N}:=(N\pi)^{s}\tilde z_{1}^{1/2}$
			at the three mesh levels $N=20,40,60$ with $n=N\nu+\mu+1$, together
			with the lower bounds
			\[
			r_1:=\Bigl(\frac{\Ab{k+1}}{\Ab{k+1+s}}\Bigr)^{1/2}
			\text{ of Lemma~\ref{lem:crude}},\qquad
			r_2:=\Bigl(\frac{\Ab{k+1-s}}{\Ab{k+1}}\Bigr)^{1/2}
			\text{ of Lemma~\ref{thm:twochan}},
			\]
			and the ratio $\theta_{60}/r_2$.
			The two bounds correspond to the two choices of test functions in the
			dual norm \eqref{eq:dualmnorm}: the single function $F_{0}^{(k+1)}$
			and the cosine sequence. Their ratio
			\begin{equation*}
				D:=\frac{r_1}{r_2}
				=\frac{\Ab{k+1}}{\bigl(\Ab{k+1-s}\,\Ab{k+1+s}\bigr)^{1/2}}\le1
			\end{equation*}
			measures how much the cosine choice improves on the single function.
			$D=1$ means that the two bounds coincide, and the smaller $D$ is, the
			larger the improvement.

			Table~\ref{tab:2} shows the general trend for $\nu\leq k$. In every
			row the ratio $\theta_{N}/r_2$ increases with $N$. In the rows with
			$D\ge0.577$ it remains below $1$ and climbs toward $1$; in the rows
			with $D\le0.449$ it reaches and exceeds $1$ and grows with $N$.
			The second
			group shows that the computed normalized error strictly exceeds $r_2$
			there, so the cosine-sequence bound is not attained. In the first
			group, and in particular in the rows with $D=1$, the computed values
			approach $r_2$ from below as $N$ grows, which is numerical evidence
			for the sharpness of the bounds.
			Table~\ref{tab:3} isolates the family $\nu=k+1$. For this family the
			lower bound \eqref{eq:nusharp} is far from sharp except for the rows with $D=1$, and
			$\theta_{60}/r_2$ tends to grow as $D$ decreases.
			Table~\ref{tab:4} makes the dependence on $D$ explicit at fixed
			$(k,s)=(6,6)$. As $\nu$ increases from $1$ to $6$, $D$ decreases
			monotonically and the ratio increases monotonically.

			These observations can be summarized as follows. When $D$ is close to
			$1$, the two choices of test functions produce nearly the same bound, and the computed
			values converge to it. This suggests that the choice of test functions is immaterial in this regime and that the bound is sharp. When $D$ is small,
			the two choices differ substantially and the computed values exceed
			the cosine bound, so even the cosine sequence does not yield a sharp bound. This is
			consistent with the observation after \eqref{eq:cosnorm} that for
			$s\ge2$ the cosine span is not dense in $\Hsp$. The supremum in the
			dual norm may be attained by a function outside its closure, and a
			better test family may exist. The last column of the tables quantifies the remaining relative gap.
		\end{itemize}

	\section{Conclusion}

	We have shown that, on a uniform one-dimensional mesh, the
	$B$-spline space is asymptotically optimal for approximating the unit ball
	of $\Hsp$ in the $L^{2}$ norm, when compared with the Kolmogorov $n$-width.
	For $k\ge s-1$ and spline spaces of lower smoothness, our lower bound for the sup-inf is asymptotically strictly larger than the $n$-width, except when $k>s-1$ and $\mu=k-2$ hold simultaneously. The numerical experiments corroborate
	these results.

	The analysis rests on two main ingredients. The first is a dual form of the
	sup-inf. It is the largest quotient of an $\Hm$-norm by an $L^{2}$-norm over
	the $L^{2}$-orthogonal complement of the spline space. The second is an
	explicit description of that complement. The $(k+1)$-st antiderivative of a function in that complement
	vanishes to order $k-\mu$ at the interior knots and to order $k+1$ at the
	endpoints. This vanishing forces the antiderivatives to oscillate at a characteristic
	frequency that decreases monotonically as the smoothness increases. In the maximally smooth case, where this frequency is smallest,  the sup-inf is asymptotically equal to the $n$-width.

	The sup-inf does not require the mesh to be uniform. Uniformity is
	what makes the analysis sharp and quantitative. The same framework can be
	adapted to measuring the error in the $H^{l}$ norm with $l\le s$, and to
	tensor-product meshes, and both are natural directions for further study.
	Finally, we note that our analysis suggests that the orthogonal complements of spaces of different smoothness may have different characteristic frequencies, which is an interesting finding that merits a deeper investigation.


\end{document}